\documentclass[11pt,a4paper,twoside]{scrartcl}
\usepackage{a4wide}
\usepackage{amsfonts}
\usepackage{amsmath}
\usepackage{amssymb}
\usepackage[numbers]{natbib}
\usepackage{booktabs}
\usepackage{algorithm}
\usepackage{algorithmic}
\usepackage{pgfplots}
\pgfplotsset{compat=newest}
\usepackage{pgfplotstable}
\usepackage{graphicx}
\usepackage{color}
\usepackage{colortbl}
\usepackage{multirow}
\usepackage{boxedminipage}
\usepackage{enumerate}
\usepackage[caption=false,font=footnotesize]{subfig}
\usepackage{rotating}
\usepackage{hyperref}
\usetikzlibrary{intersections}

\newcommand{\rn}[1]{\ensuremath{\text{\uppercase\expandafter{\romannumeral#1}}}}
\newcommand{\N}{\ensuremath{\mathbb{N}}}
\newcommand{\T}{\ensuremath{\mathbb{T}}}
\newcommand{\Z}{\ensuremath{\mathbb{Z}}}
\newcommand{\R}{\ensuremath{\mathbb{R}}}
\newcommand{\C}{\ensuremath{\mathbb{C}}}
\newcommand{\ii}{\textnormal{i}}
\newcommand{\e}{\textnormal{e}}
\newcommand{\zb}[1]{\ensuremath{\boldsymbol{#1}}}
\newcommand{\bigtimes}{\mathop{\text{\Large{$\times$}}}}

\newcommand{\rol}{{\ensuremath{\operatorname{r1l}}}}
\newcommand{\mrol}{{\ensuremath{\operatorname{mr1l}}}}
\newcommand{\boldJ}{{\ensuremath{\boldsymbol{J}}}}
\newcommand{\boldk}{{\ensuremath{\boldsymbol{k}}}}
\newcommand{\boldl}{{\ensuremath{\boldsymbol{l}}}}
\newcommand{\boldh}{{\ensuremath{\boldsymbol{h}}}}
\newcommand{\boldx}{{\ensuremath{\boldsymbol{x}}}}
\newcommand{\boldy}{{\ensuremath{\boldsymbol{y}}}}
\newcommand{\boldz}{{\ensuremath{\boldsymbol{z}}}}
\newcommand{\boldxi}{{\ensuremath{\boldsymbol{\xi}}}}
\newcommand{\boldzero}{{\ensuremath{\boldsymbol{0}}}}
\newcommand{\boldA}{{\ensuremath{\boldsymbol{A}}}}
\newcommand{\Err}{\ensuremath{\operatorname{Err}}}
\newcommand{\Res}{\ensuremath{\operatorname{Res}}}
\newcommand{\einschraenkung}{\,\rule[-5pt]{0.4pt}{14pt}\,{}}

\newtheorem{example}{Example}[section]
\newenvironment{Example}{\goodbreak \begin{example}\normalfont \rmfamily}{\bend\end{example}}

\newcommand{\bend}{\hspace*{0ex} \hfill \hbox{\vrule height
    1.5ex\vbox{\hrule width 1.4ex \vskip 1.4ex\hrule  width 1.4ex}\vrule
    height 1.5ex}}

\title{A sparse FFT approach for ODE with random coefficients}

\date{\today}
\author{Maximilian Bochmann\footnotemark[1] \and Lutz K\"ammerer\footnotemark[2] \and Daniel Potts\footnotemark[3]}

\hypersetup{pdfauthor={Maximilian Bochmann, Lutz K\"ammerer, Daniel Potts},
          pdftitle={A sparse FFT approach for ODE with random coefficients},
          pdfsubject={},
          pdfcreator = {pdflatex},
          plainpages=false,
          pdfstartview=FitH,       
          pdfview=FitH,            
          pdfpagemode=UseOutlines, 
          bookmarksnumbered=true, 
          bookmarksopen=false,     
          bookmarksopenlevel=0,   
          colorlinks=true,       
          linkcolor=black,
          citecolor=black,
          urlcolor=black}

\begin{document}

\maketitle
\begin{abstract}
\small
The paper presents a general strategy to solve ordinary differential equations (ODE), where some coefficient depend on the spatial variable and on additional random variables. The approach is based on the application of a recently developed dimension-incremental sparse fast Fourier transform. Since such algorithms require periodic signals, we discuss periodization strategies and associated necessary deperiodization modifications within the occuring solution steps.

The computed approximate solutions of the ODE depend on the spatial variable and on the random variables as well.
Certainly, one of the crucial challenges of the high dimensional approximation process is to rate the influence of each variable on the solution as well as the determination of the relations and couplings within the set of variables.
The suggested approach meets these challenges in a full automatic manner with reasonable computational costs, i.e.,
in contrast to already existing approaches, one does not need to seriously restrict the used set of ansatz functions in advance.
\medskip

\smallskip

\end{abstract}
\footnotetext[1]{ 
Chemnitz University of Technology, Faculty of Mathematics, 09107 Chemnitz, Germany,\\
  maximilian.bochmann@s2014.tu-chemnitz.de}
  
\footnotetext[2]{
  Chemnitz University of Technology, Faculty of Mathematics, 09107 Chemnitz, Germany,\\
  kaemmerer@mathematik.tu-chemnitz.de, Phone:+49-371-531-37728, Fax:+49-371-531-837728
}
\footnotetext[3]{
  Chemnitz University of Technology, Faculty of Mathematics, 09107 Chemnitz, Germany,\\
  potts@mathematik.tu-chemnitz.de, Phone:+49-371-531-32150, Fax:+49-371-531-832150
}

\section{Introduction}

During the last years, the concept of random variables has become a very popular tool to model
uncertain properties mathematically. For instance, diffusion characteristics of inhomogeneous materials
can be distinctly more accurately described by functions that additionally depend on random variables.
One common application area of these mathematical designs are diffusion coefficients in differential
equations. 
Certainly, the additional random variables affect the solvability and -- if exist -- the solutions of the differential equations under consideration. 
Besides investigations on existence, uniqueness and regularity of solutions for specific mathematical problems that involve randomness, cf. e.g. \cite{CoDeSchw10a, HaSchw13,BaCoGi17,BaCoDVGi17}, numerical solution approaches need to be developed in order to compute approximations of the desired solutions.
Accordingly, the established numerical solution approaches for differential equations without random coefficients need to be
-- at least -- extended
in order to meet the new challenges that are caused by the randomness of the diffusion coefficient.
Commonly, discretizations of the domain of the stochastic variables lead to discretized solutions that are used to compute solutions in polynomial spaces or finite element representations, cf. e.g. \cite{MaKn10,GrKuNuScheSl11,RaSchw16,HaNoTaTe16,EiPfSchn15}. One essential task in this approach is the choice of suitable polynomial spaces and corresponding basis polynomials, which can be extremely challenging for higher numbers of random variables and occuring dependencies within the random variables.
Furthermore, preferable choices of the used basis functions can improve the efficiency of the arising computations.

In the recent literature, many solution approaches deal with differential equations with random parameters in their coefficients.
Most commonly, one suggests to choose several fixed instances of the random variables and applies known solvers for the considered differential equations
without random coefficients for each of those instances. One achieves a set of solutions and computes the quantities of interest, which may be the coefficients of a specific expansion of the full solution of the differential equation or simply the expectation function, from these solutions
by applying stochastic estimators.

In this paper, we present a closed approach that deals with the spatial variable as well as the random variables simultaneously in order to solve an ordinary differential equation with a diffusion coefficient affected by randomness. In more detail, we consider the differential equation
\begin{equation}
-\frac{\partial}{\partial \eta}\left(a(\eta,\boldxi)\frac{\partial}{\partial \eta}u(\eta,\boldxi)\right)=f(\eta),\qquad u(\alpha,\boldxi)=u(\beta,\boldxi)=0 \label{eq:ODE}
\end{equation}
with homogeneous boundary conditions, where $a\colon D_a\to\R$, $D_a:=\bigtimes_{j=1}^{1+d_\boldxi}[\alpha_j,\beta_j]\subset\R\times\R^{d_\boldxi}$, is the diffusion coefficient, that depends on the spatial variable $\eta$ as well as the random variables that are the components of the vector $\boldxi\in\R^{d_\boldxi}$, and the right hand side $f\colon D_f\to\R$, $D_f:=[\alpha_1,\beta_1]\subset\R$, is a function that depends only on the spatial variable  $\eta$. Here we would like to point out, that the suggested solution approach is not restricted to the homogeneous boundary conditions or to right hand sides $f$ that are independent of the random variables $\boldxi$. Simple modifications of the presented approach lead to a solution strategy even for more general settings. However, the restrictions will simplify the notations and help to preserve clarity.

The essential restriction 
\begin{equation}
0<r\le a(\eta,\boldxi)\le R<\infty\label{eq:req_a}
\end{equation}
for all $(\eta,\boldxi)\in D_a$ guarantees the existence of a unique solution $u(\circ,\boldxi)$ of \eqref{eq:ODE} for each fixed $\boldxi$.
Hence, we suggest to approximately compute the unique solution $u$ of \eqref{eq:ODE} by means of a dimension-incremental sparse fast Fourier transform (FFT) approach,
cf. \cite{PoVo14, KaPoVo17}, and based on a direct reversion of the occurring derivatives.

On the one hand, the assumptions on the differential equation \eqref{eq:ODE} do not guarantee for periodic signals that has to be treated.
On the other hand, the sparse FFT approaches consider the input signals as periodic signals and is more successful -- in the sense of approximation rates, number of needed samples to ensure a specific accuracy, etc. --  when dealing with smooth periodic signals.
Therefore, the periodization of the arising signals will be necessary in order to compute good approximate solutions of \eqref{eq:ODE}.
At this point we would like to highlight that we approximately compute a complete solution of \eqref{eq:ODE}, which can be used to subsequently approximate several quantities of interest.
The crucial advantage of this approach is that these complete approximate solutions reveals detailed characteristics
of the random variables, i.e., the influence of each single variable on the solution as well as the interaction
between different variables. We stress that the suggested strategy automatically detects these detailed characteristics with reasonable computational costs.

The paper is organized as follows:
First we roughly outline the concept of the dimension-in\-cre\-men\-tal sparse FFT approach and indicate
the basic properties of suitable periodization mappings in Section~\ref{sec:Pre}.
Section~\ref{sec:main_results} presents the suggested strategy to treat the considered problem which leads to an approximation of
the solution of the considered differential equation.
As mentioned above, one may be interested in specific quantities of interest of this solution.
Thus, we demonstrate how to compute $n$th moments of the computed solution in Section~\ref{sec:comp_moments}.
Section~\ref{sec:numerics} contains various numerical examples, shows the operability of the suggested approach, and discusses advantages and disadvantages of the
applied sparse FFT approaches.

\section{Prerequisites}\label{sec:Pre}

\subsection{Sparse FFT}\label{Sec:sFFT}

As mentioned above, we suggest to approximate the solution of an ordinary differential equation with random coefficients using the dimension-incremental FFT approach presented in \cite{KaPoVo17}. In this section, we
declare the necessary notation and indicate the basic idea of this algorithm.

The aim of the dimension-incremental approach is the reconstruction of the Fourier coefficients $\hat{p}_\boldk$, $\boldk\in I$,
of an arbitrarily chosen trigonometric polynomial
\begin{equation}\label{poly}
p(\boldx):=\sum_{\boldk\in I} \hat{p}_\boldk \, \mathrm{e}^{2\pi\mathrm{i}\boldk\cdot\boldx},
\end{equation}
where the frequencies $\zb k$ are supported on a frequency set $I\subset\Z^d$ of finite cardinality, i.e. $|I|<\infty$.
In contrast to usual FFT algorithms, the challenge of sparse FFT algorithms is the efficient determination of the unknown frequency set $I$ in addition to the Fourier coefficients $\hat{p}_\boldk$ using only sampling values of $p$.

Appropriate thresholding strategies within dimension-incremental sparse FFT algorithms allow for the treatment of general functions $f\in L_1(\T^d)\cap\mathcal{C}(\T^d)$, i.e., the sparse FFT determines an approximation of the frequency set $I$ as well as an approximation of the (roughly) largest Fourier coefficients of the function $f$. Accordingly, the algorithms can be used in order to compute approximations
\begin{equation}
\tilde{S}_I[f](\boldx):=\sum_{\boldk\in I}\hat{f}_\boldk\e^{2\pi\ii\boldk\cdot\boldx}
\label{eq:appr_SIf}
\end{equation}
of the Fourier partial sum
$$
S_I[f](\boldx):=\sum_{\boldk\in I}c_\boldk(f)\e^{2\pi\ii\boldk\cdot\boldx}
$$
for sufficiently smooth functions $f$. In this context, the Fourier partial sum 
$S_I[f]$ is the truncated Fourier series of $f$, which implies the formal definition
of the Fourier coefficients
$$
c_\boldk(f):=\int_{\T^d}f(\boldx)\e^{-2\pi\ii\boldk\cdot\boldx}\mathrm{d}\boldx.
$$
In general, the coefficients $\hat{f}_\boldk$ are just approximations of the Fourier
coefficients $c_\boldk(f)$, since they are computed using only function evaluations of $f$ and thus are disturbed at least by aliasing.

In order to compute both, the set $I$ of the most significant frequencies $\boldk$ as well as approximations $\hat{f}_\boldk$ of the corresponding Fourier coefficients, a dimension-incremental approach was developed in \cite{PoVo14,KaPoVo17}, where the fundamental concept arises from a dimension-incremental method for the reconstruction of anharmonic trigonometric polynomials based on Prony's method, cf. \cite{PoTa12}.
An outline of this concept can be found in \cite[Sec. 2.2]{KaPoVo17}.

\begin{algorithm}[tp]
\caption{Reconstruction of a multivariate function $f$ from sampling values along (multiple) rank-1 lattices (sFFT).}\label{alg:dim_incr_sfft}
  \begin{tabular}{p{1.35cm}p{2.5cm}p{10.45cm}}
    Input:  
                & $\Gamma\subset\Z^d$ \hfill & search space in frequency domain\\
                & $f(\circ)\colon\T^d\to\C$ & function $f$ as black box (function handle) \\
                & $\theta\in\R^+$ & relative threshold \\ 
                & $s,s_\mathrm{local}\in\N$ & sparsity parameters ($s_\mathrm{local}:=s$ by default) \\ 
                & $r\in\N$ & number of detection iterations\\
                & $b\in\N$ & maximal number of multiple rank-1 lattice searches per dimension-incremental step\\
  \midrule
  \end{tabular}
  \begin{algorithmic}
	\STATE \sffamily\colorbox{black}{\bfseries\textcolor{white}{\parbox[c]{2.75cm}{\centering BLACK BOX\\ ALGORITHM}}} Available at \cite{Vo_FFTr1l}, details  in \cite[Alg. 1]{PoVo14}.
  \end{algorithmic}
  \begin{tabular}{p{1.35cm}p{2.5cm}p{10.45cm}}
  \midrule
    Output: & $I\subset\Gamma\subset\Z^d$ & set of detected frequencies, $|I|\leq \min \{s,|\Gamma|\}$ \\
            & $\boldsymbol{{\hat{f}}}\in\C^{\vert I\vert}$ & corresponding Fourier coefficients of $\tilde{S}_If$, cf. \eqref{eq:appr_SIf}
  \end{tabular}
\end{algorithm}
In this paper, we restrict the discussion to the in- and output of the algorithm, cf. Algorithm \ref{alg:dim_incr_sfft}.
We require a restricted search space $\Gamma\subset\Z^d$ in frequency domain, where the significant Fourier coefficients are assumed to be supported. For simplicity and without crucial influence on the runtime of the algorithm, we can choose a tensor product box of equal edge lengths, i.e., we fix $\Gamma=[-N,N]^d$ for a suitable edge length $2N+1$, $N\in\N$.
Since the used sampling nodes are chosen adaptively, we assume the function $f$ being given as a black box. The parameter $\theta\in\R^+$ is a thresholding for the minimal absolute values that should be accepted as significant Fourier coefficient $\hat{f}_\boldk$ and its projections in lower dimensions. Additional sparsity parameters $s,s_\mathrm{local}\in\N$
restrict the algorithm to deal with at most $s$ or $r\,s_\mathrm{local}$ frequencies in each dimension-incremental step. Here, the parameter $r$ is the number of projections that are used in each dimension-incremental step. Multiple projections are necessary in order to avoid detection failures caused by cancellations.
Since we will use only function evaluations of the function $f$ in order to compute an approximation, we have to apply suitable sampling strategies. For the case where we use multiple rank-1 lattices, the adaptive construction of the sampling set is affected by a certain small default probability. Therefore, it may happen that one has to start this construction of the sampling set more than once. The parameter $b\in\N$ can be used in order to restrict the number of restarts of the construction in each dimension-incremental step in order to guarantee the termination of the algorithm, cf. \cite{KaPoVo17}. However, this parameter 
is not restrictive during the computation, since even the choice $b=5$ is not reached in practice.

The output of Algorithm \ref{alg:dim_incr_sfft} is the frequency set $I$ and the corresponding approximated Fourier coefficients $\hat{f}_\boldk$, $\boldk\in I$, where $\tilde{S}_If$, cf. \eqref{eq:appr_SIf}, is a good approximation of $f$ when all significant frequencies are collected in $I$.

One crucial point of the dimension-incremental approach is the construction of spatial discretizations for trigonometric polynomials with frequencies in a certain, adaptively determined candidate set. Additional preferable properties of these spatial discretizations are
\begin{itemize}
\item fast discrete Fourier transform algorithms,
\item fast construction methods for the spatial discretizations, and
\item low oversampling factors, i.e., the ratio of the number of sampling values to the cardinality of the candidate set should be low.
\end{itemize}
For high dimensional sparse trigonometric polynomials the concept of multiple rank-1 lattices, cf. \cite{Kae16, Kae17} combines all these advantages, that are particularly beneficial to our targeted application.

For the sake of completeness, we give further details on the used sampling schemes.
For a given generating vector $\boldz\in\Z^d$ and a lattice size $M\in\N$, we define the rank-1 lattice
$$
\Lambda(\zb z,M)=\left\{\frac{j}{M}\boldz\bmod{\bold1}\colon j=0,\ldots,M-1\right\},
$$
where the modulo operation is applied componentwise.
For a given frequency set $I\subset\Z^d$, $|I|<\infty$, the corresponding Fourier matrix is given by
$$
\boldA=\boldA(\Lambda(\boldz,M),I):=\left(\e^{2\pi\ii\boldk\cdot\boldz\frac{j}{M}}\right)_{j=0,\ldots,M-1, \boldk\in I}.
$$
The dimension-incremental sparse FFT deals with different candidate sets of frequencies and asks for spatial discretizations for trigonometric polynomials with frequencies supported on these frequency sets. Additional requirements on $\Lambda(\boldz,M)$
guarantees the spatial discretization property, i.e., the full column rank of the matrix $\boldA(\Lambda(\boldz,M),I)$.
Due to the structure of $\Lambda(\boldz,M)$, the computations of the matrix vector products involving $\boldA$ and its pseudo inverse
can be performed by fast Fourier transform algorithms, cf. \cite{kaemmererdiss}.
These fast algorithms as well as the component--by--component construction algorithms for the used spatial discretizations, cf. \cite{Kae2013}, are the essential building blocks for the dimension-incremental sparse FFT based on single rank-1 lattices as spatial discretizations, which we denote by R1LsFFT, cf. \cite{PoVo14} for details on that approach.

A very similar approach is considered in \cite{KaPoVo17}, where the authors replaced the used sampling schemes by multiple rank-1 lattices, i.e., the spatial discretizations are constructed by the union of more than one rank-1 lattice, which provides -- at least with high probability -- asymptotically lower oversampling factors as well as much faster construction approaches for spatial discretizations. Furthermore, fast Fourier transform algorithms for the evaluation and the reconstruction of trigonometric polynomials were developed, cf. \cite{Kae17}. The corresponding dimension-incremental sparse FFT that uses these algorithms, i.e., the FFT algorithms as well as the construction algorithms for the spatial discretizations, is denoted by MR1LsFFT in the following.

Recently, a very similar
dimension-incremental sparse FFT based on random sampling was introduced in~\cite{ChIwKr18}.
One might also use this strategy in order to compute the FFT parts of the approach presented in this paper.
However, the corresponding algorithm suffers from unreasonable computational costs due to the application of direct matrix vector multiplications. For that reason, we will not use 
dimension-incremental sparse FFTs based on random sampling in our numerical tests.

The aforementioned dimension-incremental sparse FFT algorithms can be applied to periodic functions. Higher order smoothness of the treated functions often leads to smaller and thus preferable frequency sets $I$. Hence, we consider reasonable approaches to (smoothly) periodize non-periodic functions.

\subsection{Periodization}\label{ssec:periodization}

The goal of a periodization is the approximation of a non-periodic function $f\colon[\alpha,\beta]\to\C$ using
trigonometric polynomials that are naturally periodic and corresponding fast Fourier transform
algorithms.
Accordingly, we transform $f$ to a periodic function $\tilde{g}$ using a variable transform
that has the following features
\begin{align}
\varphi&\colon [0,1]\mapsto[\alpha,\beta],\nonumber\\
\varphi(0)&=\alpha,\nonumber\\
\varphi(1/2)&=\beta,\nonumber\\
\varphi(1/2-x)&=\varphi(1/2+x)\quad\textnormal{for }x\in[0,1/2],\nonumber\\
\varphi  & \textnormal{ is continuous in [0,1/2] and strictly increasing in (0,1/2)},\label{eq:basic_assumptions_periodization}
\end{align}
i.e., $f(\varphi(x))=\tilde{g}(x)$.
In more detail, we are interested in approximations of the antiderivative of the function $f$,
which leads to
\begin{align}
F(t)=\int_\alpha^t f(\tau)\mathrm{d}\tau=\int_{\varphi^{-1}(\alpha)}^{\varphi^{-1}(t)}\tilde{g}(x)\varphi'(x)\mathrm{d}x=
\int_{0}^{\varphi^{-1}(t)}\tilde{g}(x)\varphi'(x)\mathrm{d}x,\qquad t\in[\alpha,\beta].
\label{eq:integral_transform_periodization}
\end{align}
In order to compute $F(t)$, $t\in[\alpha,\beta]$, we are interested in suitable approximations
of $\tilde{g}(x)\varphi'(x)$ for $x\in[0,1/2]$, which we want to realize using trigonometric polynomials.
There are two different approaches to realize the computation of $F$.
One point of view is to approximate $\tilde{g}(x)$ and assume that $\varphi'$ is constant almost everywhere in $[0,1/2]$,
which leads to the well known tent transform approach \cite{CoKuNuSu16, SuNuCo14}.
 A more general approach will require additional assumptions
on $\varphi$ in order to obtain periodic smoothness of $\tilde{g}\varphi'$ that allow for suitable
periodic approximations. For our purposes it is enough to deal with periodizations in one dimension. In higher dimensional settings, i.e., periodizations applied to a vector of variables, we simply apply the one-dimensional periodizations to each component of the vector.

\subsubsection{Tent transform}\label{sec:tent_transform}
The so-called tent-transform \cite{CoKuNuSu16, SuNuCo14} is often used for periodization due to
its simplicity. From a geometric point of view, the tent transform
appends a mirror of the non-periodic function to the original function 
and dilates the resulting function such that its support is of length one.
In addition, the new function is shifted such that its support is exactly $[0,1]$.
In formula, the mapping
\[
\varphi\colon[0,1]\to[\alpha,\beta],\qquad\varphi(x)=\beta-|2(\beta-\alpha)(1/2-x)|
\]
realizes this periodization of a function $f\colon[\alpha,\beta]\to\C$, cf. Figure \ref{absolute} for a plot of $\varphi$, where $[\alpha,\beta]=[-1,1]$.
Certainly, this mapping $\varphi$ is not continuously differentiable.
Nevertheless, the constant first derivative $\varphi'$ within $(0,1/2)$
and $(1/2,1)$ provides advantages within the integrals that we would like to deal with.
For $t\in[\alpha,\beta]$ and $\varphi^{-1}\colon[\alpha,\beta]\to[0,1/2]$, $\varphi^{-1}(t)=\frac{t-\alpha}{2(\beta-\alpha)}$, we obtain
\begin{align*}
\int_{\alpha}^{t}f(\tau)\mathrm{d}\tau&=\int_{0}^{\varphi^ {-1}(t)}f(\varphi(x))\varphi'(x)\mathrm{d}x
=2(\beta-\alpha)\int_{0}^{\frac{t-\alpha}{2(\beta-\alpha)}}f(\varphi(x))\mathrm{d}x.
\end{align*}
Consequently, we only need to find an approximation of the antiderivative of the periodic function $f\circ\varphi$
in order to achieve an approximation of an antiderivative of the non-periodic function $f$.

\subsubsection{More general periodizations}

In addition to the basic assumptions on the peridization mapping, cf.  \eqref{eq:basic_assumptions_periodization}, we may assume periodic differentiability in order to obtain
smoother integrands in \eqref{eq:integral_transform_periodization}. Higher order smoothness of the periodization could have positive effects for the approximation of the integrand $\tilde{g}\varphi'$ using trigonometric polynomials. Roughly speaking, the smoother the function, the faster the decay of the Fourier coefficients, i.e, the smaller the cardinality of the frequency set of suitable approximating trigonometric polynomials.
In some cases it may be enough to construct periodizations of a specific fixed smoothness, since the function $f$ does not allow for higher order smoothness of $\tilde{g}$. E.g., splines of higher order seems to be ideally suited in order to guarantee the desired properties, cf. Example \ref{ex:spline}. In cases of functions $f$ of higher but unknown smoothness, infinitely differentiable mappings $\varphi$ may be an option to ensure that the periodization does not cause lower order smoothness of the integrand $\tilde{g}\varphi'$. One suitable option for such a mapping is given in Example \ref{ex:cosine}.

However, the usage of more complicated mappings may imply disadvantages in the computation of the inverse mapping of the periodization.
\begin{figure}[tb]
\subfloat[tent transform\label{absolute}]{
\begin{tikzpicture}
	\begin{axis}[
		xlabel=$x$,
		ylabel={$\varphi(x)$},
		width=.3\textwidth,
		xmax=1, xmin=0, ymax=1, ymin=-1,
		domain=0:1,
		xtick={0,0.5,1}, ytick={-1,0,1},
		yticklabel style={font=\tiny}, xticklabel style={font=\tiny},
	]
	\addplot[mark=none] {1-abs(2-4*x)};
	\end{axis}
\end{tikzpicture}
}
\subfloat[spline\label{fig:spline}]{
\begin{tikzpicture}
	\begin{axis}[
		xlabel=$x$,
		ylabel={$\varphi(x)$},
		width=.3\textwidth,
		xmax=1, xmin=0, ymax=1, ymin=-1,
		xtick={0,0.5,1}, ytick={-1,0,1},
		yticklabel style={font=\tiny}, xticklabel style={font=\tiny},
	]
	\addplot[mark=none, domain=0:0.5, samples=51] {-32*x^3+24*x^2-1};
	\addplot[mark=none, domain=0.5:1, samples=51] {32*x^3-72*x^2+48*x-9};
	\end{axis}
\end{tikzpicture}
}
\subfloat[cosine\label{fig:cosine}]{
\begin{tikzpicture}
	\begin{axis}[
		xlabel=$x$,
		ylabel={$\varphi(x)$},
		width=.3\textwidth,
		xmax=1, xmin=0, ymax=1, ymin=-1,
		domain=0:1,
		xtick={0,0.5,1}, ytick={-1,0,1},
		yticklabel style={font=\tiny}, xticklabel style={font=\tiny},
	]
	\addplot[mark=none, samples=101] {-cos(2*pi*deg(x))};
	\end{axis}
\end{tikzpicture}
}
\caption{Different periodization mappings $\varphi$.}
\end{figure}
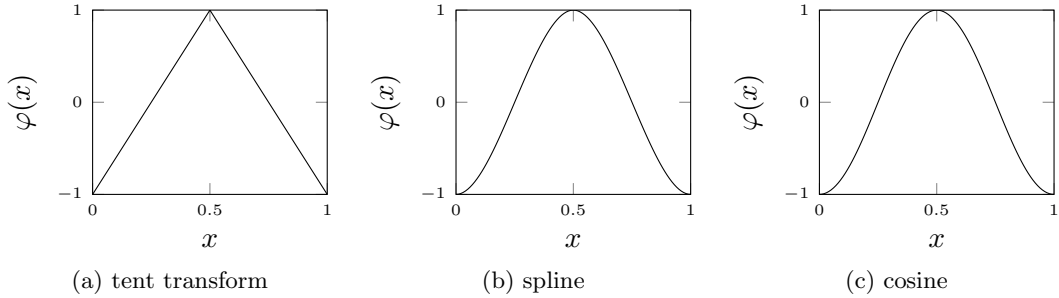
\begin{Example}\label{ex:spline}
A spline of order four can be used to construct a periodization that is two times continuously differentiable. The mapping $\varphi$, plotted in Figure \ref{fig:spline} for $[\alpha,\beta]=[-1,1]$, is given by
\begin{align*}
\varphi&\colon[0,1]\to[\alpha,\beta],\\
\varphi(x)&=
\begin{cases}
-16(\beta-\alpha)x^3+12(\beta-\alpha)x^2+\alpha& 0\le x\le 1/2,\\
16(\beta-\alpha)x^3-36(\beta-\alpha)x^2+24(\beta-\alpha)x+5\alpha-4\beta & 1/2<x\le 1.
\end{cases}
\end{align*}
\end{Example}

\begin{Example}\label{ex:cosine}
The cosine function can be used to construct an infinitely differentiable periodization mapping
\[
\varphi\colon[0,1]\to[\alpha,\beta],\qquad\varphi(x)=\frac{\alpha-\beta}{2}\cos(2\pi x)+\frac{\alpha+\beta}{2}.
\]
A corresponding plot for $[\alpha,\beta]=[-1,1]$ can be found in Figure \ref{fig:cosine}.
\end{Example}

\section{ODE solver}\label{sec:main_results}

Since the differentiations within the ODE \eqref{eq:ODE} acts on only one variable, we revert the differentiation by integration and thus obtain a formal solution
\begin{equation}
u^*(t,\boldxi)=\int_{\alpha_1}^t\frac{-\int_{\alpha_1}^\eta f(\tau)\mathrm{d}\tau+c_1(\boldxi)}{a(\eta,\boldxi)}\mathrm{d}\eta+c_2(\boldxi).
\label{eq:formal_solution1}
\end{equation}
However, for high-dimensional variables $\boldxi$ the computation of such a solution is a particular challenge.

We denote by $F$ the antiderivative of $f$ and obtain $c_2(\boldxi)=0$ since $u^*(\alpha_1,\boldxi)=0$ for homogeneous boundary conditions. The solution $u^*(t,\boldxi)$ changes to
\begin{equation}
u^*(t,\boldxi)=\underbrace{\int_{\alpha_1}^t\frac{F(\alpha_1)-F(\eta)}{a(\eta,\boldxi)}\mathrm{d}\eta}_{:=u_1(t,\boldxi)}+c_1(\boldxi)\underbrace{\int_{\alpha_1}^t\frac{1}{a(\eta,\boldxi)}\mathrm{d}\eta}_{:=u_2(t,\boldxi)},
\label{eq:formal_solution2}
\end{equation}
where we will use the term $c_1(\boldxi)$ in order to satisfy the boundary condition $u^*(\beta_1,\boldxi)=0$. In particular, requirement 
\eqref{eq:req_a} implies $u_2(\beta_1,\boldxi)>0$ and thus fixing
\begin{equation}
c_1(\boldxi):=-\frac{u_1(\beta_1,\boldxi)}{u_2(\beta_1,\boldxi)}\label{eq:formal_c1}
\end{equation}
yields homogeneous boundary conditions for $u^*$.
Accordingly, for given $f$ and $a$ we need to compute suitable approximations of $u_1$ and $u_2$.
To this end, we will apply a dimension-incremental sparse FFT approach as described in Section \ref{Sec:sFFT}. Since these FFT algorithms handles periodic signals, we need to periodize the upcoming functions.

\subsection{Integration of the right hand side  $f$}
First, we determine the term $F(\alpha_1)-F(t)$ in \eqref{eq:formal_solution2} from above by approximating and integrating $f$.
To this end, we periodize $f$ using a suitable periodization $\varphi$, cf. Section \ref{ssec:periodization},
\[
\tilde{f}(x)=f(\varphi(x)).
\]
Accordingly, we obtain
\[
\int f(\tau)\mathrm{d}\tau=\int f(\varphi(x))\varphi'(x)\mathrm{d}x=\int \tilde{f}(x)\varphi'(x)\mathrm{d}x
\]
and approximate the integrand on the right hand side by a trigonometric polynomial \begin{equation}
\tilde{S}_N\left[\tilde{f}\varphi'\right](x)=\sum_{k=-N}^{N}\hat{a}_k\e^{2\pi\ii kx}.\label{eq:S_N_f_of_varphi_times_varphiprime}
\end{equation}
An antiderivative of $\tilde{S}_N\left[\tilde{f}\varphi'\right]$ is given by 
\[\int\tilde{S}_N\left[\tilde{f}\varphi'\right](x)\mathrm{d}x=\hat{a}_0x+\sum_{1\le|k|\le N}\frac{\hat{a}_k}{2 k\pi\ii}\e^{2\pi\ii kx},\]
which yields 
\[
F(\tau)+c=
\int f(\tau)\mathrm{d}\tau \approx\hat{a}_0\varphi^{-1}(\tau)+\sum_{1\le|k|\le N}\frac{\hat{a}_k}{2 k\pi\ii}\e^{2\pi\ii k\varphi^{-1}(\tau)}.
\]
Consequently, we denote the approximation of the term $F(\alpha_1)-F(t)$ by $\breve{F}(t)$ and obtain
\begin {align}
F(\alpha_1)-F(\eta)\approx \breve{F}(\eta)&:=\hat{a}_0(\varphi^{-1}(\alpha_1)-\varphi^{-1}(\eta))\nonumber\\
&\qquad\qquad+\sum_{1\le|k|\le N}\frac{\hat{a}_k}{2 k\pi\ii}\e^{2\pi\ii k\varphi^{-1}(\alpha_1)}-\sum_{1\le|k|\le N}\frac{\hat{a}_k}{2 k\pi\ii}\e^{2\pi\ii k\varphi^{-1}(\eta)}\nonumber\\
&=-\hat{a}_0\varphi^{-1}(\eta)-\sum_{1\le|k|\le N}\frac{\hat{a}_k}{2 k\pi\ii}(\e^{2\pi\ii k\varphi^{-1}(\eta)}-1).\label{eq:def_breveF}
\end{align}
\subsection{Approximating $u_1$, $u_2$, and $c_1$}
\label{sssec:approx_u1_u2_c1}

We denote the integrands in \eqref{eq:formal_solution2} that determine $u_1$ and $u_2$ by $v_1$ and $v_2$, respectively, i.e. we have
\begin{align*}
v_1(\eta,\boldxi)&:=\frac{F(\alpha_1)-F(\eta)}{a(\eta,\boldxi)}
\qquad\textnormal{and}\qquad
v_2(\eta,\boldxi):=\frac{1}{a(\eta,\boldxi)}.
\end{align*}
First we consider the function $v_1$ and plug in the approximation $\breve{F}(\eta)$ of 
$F(\alpha_1)-F(\eta)$ from~\eqref{eq:def_breveF}. This yields an approximation of $v_1$
\[
\breve{v}_1(\eta,\boldxi):=\frac{\breve{F}(\eta)}{a(\eta,\boldxi)}.
\]
Now, our goal is to construct an antiderivative of $\breve{v}_1$ with respect to $\eta$, which is an approximation of the antiderivative of $v_1$. 
To this end, we construct a periodization of $\breve{v}_1$ using mappings $\varphi_\eta$ and $\varphi_\boldxi$. We take into account the influence of the periodization during integration with respect to the first variable, which leads to the periodic integrand
$\tilde{\breve{v}}_1(x,\boldy)=\breve{v}_1(\varphi_\eta(x),\varphi_\boldxi(\boldy))\,\varphi_\eta'(x)$.
We compute a corresponding approximation
\begin{equation}
\tilde{S}_{I_1}[\tilde{\breve v}_1](x,\boldy):=\sum_{(k,\boldl)\in I_1\subset\Z^{1+d_\boldy}}\hat{b}_{(k,\boldl)}\e^{2\pi\ii(kx+\boldl\cdot\boldy)},\label{eq:S_I_of_v1}
\end{equation}
which can be done by sparse FFT approaches, as described in Section \ref{Sec:sFFT}, similar to those described in \cite{PoVo14,KaPoVo17}. The antiderivative of $\tilde{S}_{I_1}[\tilde{\breve v}_1]$ with respect to $x$ is given by
\begin{align}
\int\tilde{S}_{I_1}[\tilde{\breve v}_1](x,\boldy)\mathrm{d}x=\sum_{\substack{(k,\boldl)\in I_1\\k\neq 0}}\frac{\hat{b}_{(k,\boldl)}}{2k\pi\ii}\e^{2\pi\ii(kx+\boldl\cdot\boldy)}+x\sum_{(0,\boldl)\in I_1}\hat{b}_{(0,\boldl)}\e^{2\pi\ii\boldl\cdot\boldy} + C(\boldy)\label{eq:antiderivative_tildebrevev1}
\end{align}
We choose $C(\boldy):=-\sum_{\substack{(k,\boldl)\in I_1\\k\neq 0}}\frac{\hat{b}_{(k,\boldl)}}{2k\pi\ii}\e^{2\pi\ii\boldl\cdot\boldy}$
in order to guarantee $\breve{u}_1(\alpha_1,\boldxi)=0$ for all $\boldxi\in D_a':=\bigtimes_{j=2}^{d_\boldxi+1}[\alpha_j,\beta_j]$, and we roll the periodization back, which leads to the approximation
\begin{equation}
\breve{u}_1(t,\boldxi):=\sum_{\substack{(k,\boldl)\in I_1\\k\neq 0}}\underbrace{\frac{\hat{b}_{(k,\boldl)}}{2k\pi\ii}}_{=:\mathfrak{b}_{(k,\boldl)}}\left(\e^{2\pi\ii k\varphi_\eta^{-1}(t)}-1\right)\e^{2\pi\ii\boldl\cdot\varphi_\boldxi^{-1}(\boldxi)}+\varphi_\eta^{-1}(t)\sum_{(0,\boldl)\in I_1}\underbrace{\hat{b}_{(0,\boldl)}}_{=:\mathfrak{b}_{(0,\boldl)}}\e^{2\pi\ii\boldl\cdot\varphi_\boldxi^{-1}(\boldxi)}\label{eq:breve_u1}
\end{equation}
of $u_1$ given in \eqref{eq:formal_solution2}.

The analogous approach, but without approximating $v_2$, leads to an approximation of $u_2(t,\boldxi)$
\begin{equation}
\breve{u}_2(t,\boldxi):=\sum_{\substack{(k,\boldl)\in I_2\subset \mathbb Z^{1+d_\boldy}\\k\neq 0}}\underbrace{\frac{\hat{c}_{(k,\boldl)}}{2k\pi\ii}}_{=:\mathfrak{c}_{(k,\boldl)}}\left(\e^{2\pi\ii k\varphi_\eta^{-1}(t)}-1\right)\e^{2\pi\ii\boldl\cdot\varphi_\boldxi^{-1}(\boldxi)}+\varphi_\eta^{-1}(t)\sum_{(0,\boldl)\in I_2}\underbrace{\hat{c}_{(0,\boldl)}}_{=:\mathfrak{c}_{(0,\boldl)}}\e^{2\pi\ii\boldl\cdot\varphi_\boldxi^{-1}(\boldxi)}.\label{eq:breve_u2}
\end{equation}
The construction of $\breve{u}_j(t,\boldxi)$ yields $\breve{u}_j(\alpha_1,\boldxi)=0$, $j=1,2$.
Consequently, each linear combination of $\breve{u}_1$ and $\breve{u}_2$ satisfies the homogeneous boundary condition in $t=\alpha_1$.
A suitable approximation of $c_1(\boldxi)$, cf. \eqref{eq:formal_c1}, will lead to  a linear combination of $\breve{u}_1$ and $\breve{u}_2$
that also satisfies the homogeneous boundary condition in $t=\beta_1$.
To this end, we periodize $\breve{u}_1$ as well as $\breve{u}_2$ and construct the approximation
\[
\breve{c}_1(\boldxi):=-\frac{\breve{u}_1(\beta_1,\boldxi)}{\breve{u}_2(\beta_1,\boldxi)}
\qquad\text{and its periodization}\qquad
\tilde{\breve{c}}_1(\boldy)=-\frac{\breve{u}_1(\varphi_\eta(1/2),\varphi_\boldxi(\boldy))}{\breve{u}_2(\varphi_\eta(1/2),\varphi_\boldxi(\boldy))},
\]
which are well defined due to the requirements on the diffusion coefficient $a$. We stress on the fact that the periodizations of $\breve{u}_j$
do not coincide to the terms in \eqref{eq:antiderivative_tildebrevev1}, since these are non-periodic in general due to the terms that are linear in $x$.

We approximate $\tilde{\breve c}_1$ using sparse FFT approaches by
\begin{equation}
\tilde{S}_{I_3}[\tilde{\breve c}_1](\boldy):=\sum_{\zb l\in I_3\subset\Z^{d_\boldy}}\mathfrak{d}_\boldl\e^{2\pi\ii\boldl\cdot\boldy}.\label{eq:S_I_c}
\end{equation}
and achieve an approximation of the non-periodic function $c_1$ by
\begin{equation}
\breve{\breve c}_1(\boldxi)=\sum_{\zb l\in I_3}\mathfrak{d}_\boldl\e^{2\pi\ii\boldl\cdot\varphi_\boldxi^{-1}(\boldxi)}.\label{eq:breve_breve_c1}
\end{equation}

Altogether, an approximation of $u^*(t,\boldxi)$, cf. \eqref{eq:formal_solution2}, is then given by
\begin{align}
\breve{u}(t,\boldxi):=\breve{u}_1(t,\boldxi)+\breve{\breve c}_1(\boldxi)\,\breve{u}_2(t,\boldxi),
\label{eq:breve_u}
\end{align}
which actually is built of three Fourier series combined with inverse mappings of the periodizations $\varphi_\eta$ and $\varphi_\boldxi$.
Algorithm \ref{alg:central_strategy} summarizes the approach stated above. 

\begin{algorithm}[tb]
\caption{Basic procedure for computing an approximation of the solution of \eqref{eq:ODE} using a dimension-incremental sparse FFT approach}\label{alg:central_strategy}
  \begin{tabular}{p{1.35cm}p{4.9cm}p{8.05cm}}
    Input:      & $f\colon\T\to\C$ & function handle of right hand side $f$\\
                & $a\colon\T^{1+d_\boldxi}\to\C$ & function handle of random coefficient $a$\\
                & $\varphi_\eta\colon [0,1]\to[\alpha_1,\beta_1]$ & periodization mapping of spatial variable $\eta$\\
                & $\varphi_\eta^{-1}\colon [\alpha_1,\beta_1]\to [0,1/2]$ & inverse of the periodization mapping $\varphi_\eta$\\
                & $\varphi_\eta'\colon [0,1]\to\R$ & first derivative of $\varphi_\eta$\\
                & $\varphi_\boldxi\colon [0,1]^{d_\boldxi}\to\bigtimes_{j=2}^{d_{\boldxi}+1}[\alpha_j,\beta_{j}]$ & periodization mapping of random variables $\boldxi$\\
                & $N\in\N$, $\theta\in\R$, $s\in\N$ & sFFT parameters\\
                \midrule
  \end{tabular}
  \begin{algorithmic}[1]
        \STATE Compute the Fourier coefficients $\left\{\hat{a}_k\right\}_{k=-N}^N$ of $\tilde{S}_N\left[\tilde{f}\varphi'\right]$\\\qquad by means of a 1d FFT using function values of $(f\circ\varphi_\eta)\varphi_\eta'$, cf. \eqref{eq:S_N_f_of_varphi_times_varphiprime}\label{alg:censtra_approx_rhs}
		\STATE Compute the coefficients of the finite sum representation of non-periodic $\breve{F}$\\\qquad by modifying the coefficients $\left\{\hat{a}_k\right\}_{k=-N}^N$, cf. \eqref{eq:def_breveF}
		\STATE Compute the Fourier coefficients $\{\hat{b}_{(k,\boldl)}\}_{(k,\boldl)\in I_1}$ of $\tilde{S}_{I_1}[\tilde{\breve v}_1]$ by means of an sFFT algorithm\\\qquad
		using sampling values of $\breve{F}(\varphi_\eta(x))\varphi_\eta'(x)/a(\varphi_\eta(x),\varphi_\boldxi(\boldy))$ , cf. \eqref{eq:S_I_of_v1}
		\label{alg:censtra_approx_v1}
		\STATE Compute the Fourier coefficients $\{\hat{c}_{(k,\boldl)}\}_{(k,\boldl)\in I_2}$ of $\tilde{S}_{I_2}[\tilde{\breve v}_2]$ by means of an sFFT algorithm\\\qquad
		using sampling values of $\varphi_\eta'(x)/a(\varphi_\eta(x),\varphi_\boldxi(\boldy))$ , similar to \eqref{eq:S_I_of_v1}
		\label{alg:censtra_approx_v2}
		\STATE Compute the coefficients $\{\mathfrak{b}_{(k,\boldl)}\}_{(k,\boldl)\in I_1}$ and $\{\mathfrak{c}_{(k,\boldl)}\}_{(k,\boldl)\in I_2}$\\
		\qquad of the finite sum representation of non-periodic $\breve{u}_1$ and $\breve{u}_2$\\
		\qquad by modifying the coefficients $\{\hat{b}_{(k,\boldl)}\}_{(k,\boldl)\in I_1}$ and $\{\hat{c}_{(k,\boldl)}\}_{(k,\boldl)\in I_2}$, cf. \eqref{eq:breve_u1} and \eqref{eq:breve_u2}
\STATE Compute the Fourier coefficients $\{\mathfrak{d}_{\boldl}\}_{\boldl\in I_3}$ of $\tilde{S}_{I_3}[\tilde{\breve c}_1]$ by means of an sFFT algorithm\\\qquad
		using sampling values of $\breve{u}_1(\varphi_\eta(1/2),\varphi_\boldxi(\boldy))/\breve{u}_2(\varphi_\eta(1/2),\varphi_\boldxi(\boldy))$ , cf. \eqref{eq:S_I_c}
		
\end{algorithmic}
  \begin{tabular}{p{1.35cm}p{4.9cm}p{8.05cm}}
  \midrule
    Output: & $\{\mathfrak{b}_{(k,\boldl)}\}_{(k,\boldl)\in I_1}$ & coefficients of $\breve{u}_1$, cf. \eqref{eq:breve_u1} \\
    & $\{\mathfrak{c}_{(k,\boldl)}\}_{(k,\boldl)\in I_2}$ &  coefficients of $\breve{u}_2$, cf. \eqref{eq:breve_u2} \\
    & $\{\mathfrak{d}_{\boldl}\}_{\boldl\in I_3}$ &  coefficients of $\breve{\breve{c}}_1$, cf. \eqref{eq:breve_breve_c1}
  \end{tabular}
\end{algorithm}

\subsection{Tent transform in spatial domain}

For the specific choice of the tent transform, cf. Section \ref{sec:tent_transform}, for the periodization $\varphi_\eta$ and componentwise for $\varphi_\boldxi$, we obtain some simplifications in the calculations above. Moreover, uniformly distributed random variables $\xi_j$ lead to additional simplifications due to the constant probability density. We observe

\begin{align*}
\int_{\alpha_1}^\eta f(\tau)\mathrm{d}\tau&=2\int_0^{\frac{\eta-\alpha_1}{2(\beta_1-\alpha_1)}}f\big(\beta_1-|2(\beta_1-\alpha_1)(1/2-x)|\big)\mathrm{d}x\\
\intertext{and, hence, it is enough to compute an approximation of $\tilde{f}=f\circ\varphi_\eta$, due to the equality}
\tilde{S}_N[\tilde f\varphi_\eta']\einschraenkung_{[0,1/2]}&=2\tilde{S}_N[\tilde f]\einschraenkung_{[0,1/2]}.
\end{align*}
The approximation of $\tilde{S}_N[\tilde f]$ is preferable, since $\tilde{f}\varphi_\eta'$ is not continuous in the case where $f(\beta_1)\neq 0$ and, thus, problematic to approximate using trigonometric polynomials.
Subsequent to the computation of $\tilde{S}_N[\tilde f]:=\sum_{k=-N}^N\hat{a}_k\e^{2\pi\ii k x}$, the calculations of the antiderivative and deperiodization leads to
\begin{align*}
\breve{F}(\eta)&:=2\left(-\hat{a}_0\varphi^{-1}(\eta)+\sum_{1\le|k|\le N}\frac{\hat{a}_k}{2 k\pi\ii}-\sum_{1\le|k|\le N}\frac{\hat{a}_k}{2 k\pi\ii}\e^{2\pi\ii k\varphi^{-1}(\eta)}\right)\nonumber\\
&=-\hat{a}_0\frac{\eta-\alpha_1}{\beta_1-\alpha_1}-\sum_{1\le|k|\le N}\frac{\hat{a}_k}{k\pi\ii}(\e^{\pi\ii k\frac{\eta-\alpha_1}{\beta_1-\alpha_1}}-1).
\end{align*}

The approximations $\breve{u}_1$ and $\breve{u}_2$, cf. \eqref{eq:breve_u1} and \eqref{eq:breve_u2}, can be computed in the exact same manner.
Altogether, using the tent transform in spatial domain requires slight modifications in Algorithm \ref{alg:central_strategy} in lines \ref{alg:censtra_approx_rhs}, \ref{alg:censtra_approx_v1}, and \ref{alg:censtra_approx_v2}. The used sampling values must be computed with a factor $2$ instead of $\varphi_\eta'$.

\section{Computing moments  of the solution}
\label{sec:comp_moments}

In Section \ref{sec:main_results}, we discussed a strategy for computing an approximate solution $\breve{u}$, cf. \eqref{eq:breve_u}, of the ODE in \eqref{eq:ODE}.
The computation of quantity of interests needs some further investigations.
For simplicity, we demonstrate one approach to compute approximations of the $n$th moments of the solution $u^*$ of the ODE in \eqref{eq:ODE}
based on the approximation $\breve{u}$.
To this end, we denote the domain of the random variables by
$D_a':=\bigtimes_{j=2}^{d_\boldxi+1}[\alpha_j,\beta_j]$.

The $n$th moment of the solution $u^*$ of \eqref{eq:ODE} is given by
\[
u^*_{\operatorname{E}^n}(t):=\operatorname{E}(\left(u^*(t,\circ)\right)^n)=\int_{D_a'}\left(u^*(t,\boldxi)\right)^n\mathrm{d}\mu(\boldxi)=\int_{D_a'}\left(u^*(t,\boldxi)\right)^n\rho(\boldxi)\mathrm{d}\boldxi,
\]
where $\rho$ is the probability density function of the random variable vector $\boldxi$.
Periodization yields
\begin{align}
\int_{D_a'}\left(u^*(t,\boldxi)\right)^n\rho(\boldxi)\mathrm{d}\boldxi&=
\int_{[0,1/2]^{d_{\boldxi}}}\left(u^*(\varphi_t(x),\varphi_\boldxi(\boldy))\right)^n\rho(\varphi_\boldxi(\boldy))\left|\det(\boldJ)\right|\mathrm{d}\boldy,\nonumber
\intertext{where $\boldJ$ is the involved Jacobian matrix. Assuming $\varphi_\boldxi$ is a periodization that acts on each component of $\boldxi$ separately, cf. Section \ref{ssec:periodization}, the determinant of the Jacobian matrix is a tensor product function and we continue}
&=2^{-d_\boldxi}\int_{\T^{d_\boldxi}}\left(u^*(\varphi_t(x),\varphi_\boldxi(\boldy))\right)^n\rho(\varphi_\boldxi(\boldy))\prod_{j=1}^{d_\boldxi}\left|\varphi_{\xi_j}'(y_j)\right|\mathrm{d}\boldy\nonumber\\
&\approx 2^{-d_\boldxi}\int_{\T^{d_\boldxi}}\underbrace{\left(\breve{u}(\varphi_t(x),\varphi_\boldxi(\boldy))\right)^n\rho(\varphi_\boldxi(\boldy))\prod_{j=1}^{d_\boldxi}\left|\varphi_{\xi_j}'(y_j)\right|}_{=:w_n(x,\boldy)}\mathrm{d}\boldy
.
\label{eq:nth_moment_w}
\end{align}
We approximate the integrand $w_n$ using a sparse FFT approach and achieve a Fourier partial sum
\[
S_{I_4}[w_n](x,\boldy)=\sum_{(k,\boldl)\in I_4}\hat{a}_{(k,\boldl)}\e^{2\pi\ii(kx+\boldl\cdot\boldy)}.
\]
Integrating $S_{I_4}[w_n]$ instead of $w_n$ in \eqref{eq:nth_moment_w} leads to the approximation
\begin{align*}
\breve{u}^*_{\operatorname{E}^n}(t)&:=
2^{-d_\boldxi}\int_{\T^{d_\boldxi}}\sum_{(k,\boldl)\in I_w}\hat{a}_{(k,\boldl)}\e^{2\pi\ii(kx+\boldl\cdot\boldy)}\mathrm{d}\boldy
=2^{-d_\boldxi}\sum_{(k,\boldzero)\in I_4}\hat{a}_{(k,\boldzero)}\e^{2\pi\ii kx}\\
&=2^{-d_\boldxi}\sum_{(k,\boldzero)\in I_4}\hat{a}_{(k,\boldzero)}\e^{2\pi\ii k\varphi_t^{-1}(t)}
\end{align*}
of $u^*_{\operatorname{E}^n}$
since each monomial that depends on $\boldy$ integrates to zero.

\section{Numerical results}
\label{sec:numerics}

For our numerical tests, we use an example from \cite{BoRaSchw17}.
The goal is to numerically solve the boundary-value problem

\begin{align}
- \frac{\partial}{\partial \eta}\left( a(\eta,\boldxi)\frac{\partial}{\partial \eta}u(\eta,\boldxi) \right)=10,\text{ with } u\equiv0 \text{ at } \partial(0,1),
\label{eq:num_ode}
\end{align}

where the random coefficient $a:[0,1] \times [-1,1]^{d_{\boldxi}} \rightarrow \mathbb{R}$ is given by

\[
a(\eta,\boldxi)=a_{0}+\sum_{j=1}^{d_{\boldxi}/2} \xi_{2j-1} \frac{\cos(j \pi \eta)}{j^{\gamma}}+ \xi_{2j} \frac{\sin(j \pi \eta)}{j^{\gamma}}, 
\]
with $\gamma \in \mathbb{R}$, $\gamma > 1$, $a_{0}\in \mathbb{R}$, $a_{0}>2 \zeta(\gamma)$, $d_{\boldxi} \in 2\mathbb{N}$, and $\zeta$ denotes the Riemann zeta function.
The random coefficient $a$ is bounded in the interval $\left[a_{0}- 2\zeta(\gamma),a_{0} + 2\zeta(\gamma)\right]$ and thus the differential
equation \eqref{eq:num_ode} is uniquely solvable for fixed $\boldxi \in [-1,1]^{d_{\boldxi}}$.\\ 

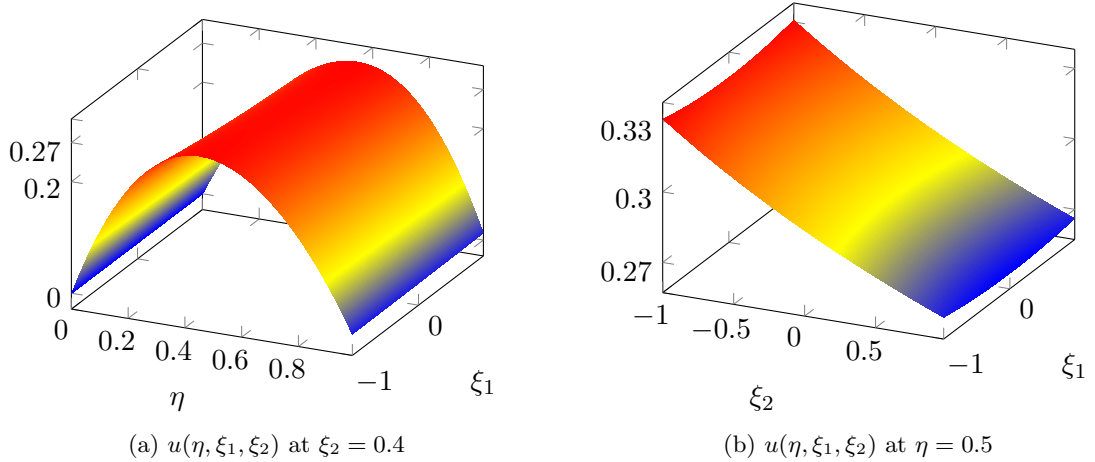
\begin{figure}
\centering
\hfill
\captionsetup[subfigure]{oneside,margin={0.5cm,0cm}}
\subfloat[$u(\eta,\xi_{1},\xi_{2})$ at $\xi_{2}=0.4$]{
\begin{tikzpicture}
\begin{axis}[
xlabel=$\eta$,
ylabel=$\xi_1$,
extra z ticks={0.27},
width=0.45\textwidth
]
 \addplot3[surf,mesh/rows=34,mesh/ordering=colwise,shader=interp] file {P1.dat};
\end{axis}
\end{tikzpicture}}
\captionsetup[subfigure]{oneside,margin={0.5cm,0cm}}
\hfill
\subfloat[$u(\eta,\xi_{1},\xi_{2})$ at $\eta=0.5$]{
\begin{tikzpicture}
\begin{axis}[
xlabel=$\xi_2$,
ylabel=$\xi_1$,
extra z ticks={0.27,0.33},
width=0.45\textwidth
]
 \addplot3[surf,mesh/rows=34,mesh/ordering=colwise,shader=interp] file {P.dat};
\end{axis}
\end{tikzpicture}}
\hfill
\caption{Solution $u(\eta,\xi_{1},\xi_{2})$ of \eqref{eq:num_ode} for $d_{\boldxi}=2$ parameters and $a_{0}=4.3$, $\gamma=2$.}\label{fig:example_plots_solution}
\end{figure}

The parameters $\xi_{k}, \, k=1,\dots,d_{\boldxi}$ can be interpreted as random variables. Here we choose them to be uniformly distributed 
\[
\xi_{k} \sim U([-1,1])
\]
and we fix $a_{0}=4.3$ and $\gamma=2$. In Figure \ref{fig:example_plots_solution} we
(partially) plotted an approximation of the solution of this differential equation,
where we restricted the number of random variables to two.

In order to demonstrate the applicability of the presented approach, we specify the settings of the applied algorithmic components.
On the one hand, we restrict the numerical tests to the tent transform as periodization mapping, cf. Section \ref{sec:tent_transform}, since this seems to be the most unfavourable choice due to its relatively low smoothness.
On the other hand, we have to specify the applied sparse FFT approaches and the corresponding parameters.
We choose three different sparsity levels $s$ and refinements $N$ for our approximated solutions $\breve{u}^{\rho}$, $\rho=\rn{1},\rn{2},\rn{3}$, cf. Section \ref{sssec:approx_u1_u2_c1}. Furthermore, we apply two different algorithms for computing approximate solutions denoted by $\breve{u}_{\rol}^{\rho}$ and $\breve{u}_{\mrol}^{\rho}$ namely the sFFT-algorithms that use sampling schemes that are \emph{r}ank-\emph{1} \emph{l}attices and \emph{m}ultiple \emph{r}ank-\emph{1} \emph{l}attices, respectively.
We call the corresponding sFFT algorithms R1LsFFT and MR1LsFFT.
The basic structure of both algorithms is described in \cite[Alg. 1]{PoVo14}. The crucial differences of the R1LsFFT and the MR1LsFFT are in step 2b and 2f, where the first approach uses the component--by--component construction as described in ``Algorithm 1'' in \cite[Sec. 2.2.1]{PoVo14} in order to determine suitable generating vectors and the latter approach uses \cite[Alg. 4]{Kae17} with $c=2$ in order to determine multiple rank-1 lattice discretizations.
\\
Tables \ref{tab:parameter_and_nosn_sFFT} and \ref{tab:parameter_and_nosn_msFFT} show the parameters we used in columns two to five for the two different sFFT algorithms. The impacts of these parameters are shortly described in Section~\ref{Sec:sFFT}.

\begin{table}[tb]
\begin{center}
\pgfplotstabletypeset
[every head row/.style={before row=\toprule,after row=\midrule},every last row/.style={after row=\bottomrule},columns/j/.style={column name=$\rho$, string type},columns/theta/.style={column name=$\theta$},columns/rtilde/.style={column name=$r$, column type=c|},columns/M_v/.style={column name=$M_{\tilde{S}_{I_1}[\tilde{\breve v}_1]}$},columns/M_w/.style={column name=$M_{\tilde{S}_{I_2}[\tilde{v}_2]}$},columns/M_u/.style={column name=$M_{S_{I_{4}}[w_1]}$,column type=c|},columns/Gamma/.style={column name=$|\hat{G}_{N}^{21}|$}
]
{j N s theta rtilde M_v M_w M_u Gamma
 {${\rn{1}}$} 32 1000 1e-12 5 27125692 29611654 7287006 117809547936177440979200839996337890625 
 {${\rn{2}}$} 64 5000 1e-12 5 566884162 715262588 172068836 210079616665422015609148165659772062692313729
 {${\rn{3}}$} 128 8000 1e-12 5 2464232808 3058109298 697944658 406065029957821485337176148033167189926689985664257
}
\end{center}
\caption{Parameter settings, the number of samples $M$ used for the computation of the approximations of $\tilde{\breve{v}}_{1}$, $\tilde{v}_{2}$, and $S_{I_4}[w_1]$ for $d_\boldxi=20$, cf. Sections \ref{sssec:approx_u1_u2_c1} and \ref{sec:comp_moments}, and the total cardinality of the box of frequency candidates for the sFFT algorithm that uses rank-1 lattices as spatial discretizations for $d_\boldxi=20$.
} 
\label{tab:parameter_and_nosn_sFFT}
\end{table}

\begin{table}[tb]
\begin{center}
\pgfplotstabletypeset
[every head row/.style={before row=\toprule,after row=\midrule},every last row/.style={after row=\bottomrule},columns/j/.style={column name=$\rho$, string type},columns/theta/.style={column name=$\theta$},columns/rtilde/.style={column name=$r$, column type=c|},columns/M_v/.style={column name=$M_{\tilde{S}_{I_1}[\tilde{\breve v}_1]}$},columns/M_w/.style={column name=$M_{\tilde{S}_{I_2}[\tilde{v}_2]}$},columns/M_u/.style={column name=$M_{S_{I_{4}}[w_1]}$, column type=c|},columns/Gamma/.style={column name=$|\hat{G}_{N}^{21}|$}
]
{j N s theta rtilde M_v M_w M_u Gamma
 {$\rn{1}$} 32 1000 1e-12 5 110061791 115366757 25463901 117809547936177440979200839996337890625 
 {$\rn{2}$} 64 5000 1e-12 5 1346337799 1446647649 353826741 210079616665422015609148165659772062692313729
 {$\rn{3}$} 128 8000 1e-12 5 4495557601 4750274729 848629827 406065029957821485337176148033167189926689985664257
}
\end{center}
\caption{Parameter settings, the number of samples $M$ used for the computation of the approximations of $\tilde{\breve{v}}_{1}$, $\tilde{v}_{2}$, and $S_{I_{4}}[w_1]$ for $d_\boldxi=20$, cf. Sections \ref{sssec:approx_u1_u2_c1} and \ref{sec:comp_moments}, and the total cardinality of the box of frequency candidates for the sFFT algorithm that uses multiple rank-1 lattices as spatial discretizations for $d_\boldxi=20$.
} 
\label{tab:parameter_and_nosn_msFFT}
\end{table}

Increasing the number $d_\boldxi$ of random variables yields approximation problems of higher dimensionality. Clearly for practical applications, the number of random variables needs to be suitably bounded. The used diffusion coefficient $a$ is build in such a way, that the influence of the random variable $\xi_j$ decreases with growing index $j$. Our first crucial task is to estimate the index $j$ for
which we can truncate the series expansion of $a$ without losing significant information of $a$. In other words we would like to estimate a suitable number $d_\boldxi$.

\begin{Example}
To this end, we computed the approximation of a solution of \eqref{eq:num_ode} by our approach with a fixed large number $d_\boldxi=40$ of random variables, i.e., we treat a 41-dimensional approximation problem. 
We end up with an approximation $\breve{u}_{\rol}^{\rn{3}}$ as represented in \eqref{eq:breve_u}.
In order to simplify the considerations on the influence and the interactions on the variables of $\breve{u}$ we apply periodizations and the sparse FFT approach on $\breve{u}$ which leads in essence to
a single Fourier sum representation of $\breve{u}$. 
The associated frequency set of this approximate solution -- together with the absolute values of the occuring (Fourier) coefficients of this solution -- allow for rating the random variables to their importance. In particular, if the expansion $h_j-l_j$, $(k,\boldh),(k',\boldl)\in I$ of the frequency set in direction $j$ is zero -- or very small and the corresponding coefficients almost zero in relation to the largest occuring coefficients -- the solution does not or not significantly depend on the variable $\xi_j$. Accordingly, leaving out this variable should not cause significant errors.
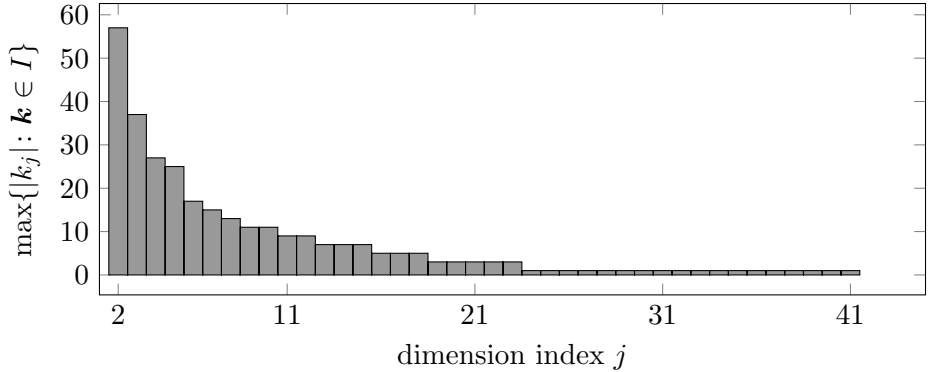
\begin{figure}[tb]
\centering
\begin{tikzpicture}
\begin{axis}[ybar stacked,xlabel={dimension index $j$},ylabel={$\max\{|k_j|\colon\zb k\in I\}$},bar width=1,xmin=0,
ytick={0,10,20,30,40,...,128},
xtick={1,10,20,30,40},
xticklabels={2,11,21,31,41},
width=.8\textwidth, height=.35\textwidth]
\addplot[color=black, fill=black!40] coordinates
	{(1,57) (2,37) (3,27) (4,25) (5,17) (6,15) (7,13) (8,11) (9,11) (10,9) (11,9) (12,7) (13,7) (14,7) (15,5) (16,5) (17,5) (18,3) (19,3) (20,3) (21,3) (22,3) (23,1) (24,1) (25,1) (26,1) (27,1) (28,1) (29,1) (30,1) (31,1) (32,1) (33,1) (34,1) (35,1) (36,1) (37,1) (38,1) (39,1) (40,1)  };
\end{axis}
\end{tikzpicture}
\caption{Directional expansion of the frequency set $I\subset\Z^{1+d_\boldxi}$ of $\breve{u}$ for the random variables.}
\label{fig:expansion_d40}
\end{figure}
Figure \ref{fig:expansion_d40} indicates the expansions in each coordinate direction of the frequency set of $\breve{u}$ for $d_\boldxi=40$. Obviously, the last 18 random variables have a very small expansion. We stress that the variables $\xi_{21}$ as well as $\xi_{22}$ have a significant frequency support but can be neglected due to the low order of magnitude of its Fourier coefficients.
For these reasons, we restrict the number of random variables to $d_\boldxi=20$ in the following
experiments.
\end{Example}

As mentioned in the last example, we fix $d_\boldxi=20$. We solved \eqref{eq:num_ode} by the means of the sparse FFT approaches that uses single or multiple rank-1 lattices as spatial discretizations. The applied parameter constellations are presented in Tables \ref{tab:parameter_and_nosn_sFFT} and \ref{tab:parameter_and_nosn_msFFT}. Both tables contains
the total amount of samples that were used for the approximation of the functions $\tilde{\breve{v}}_{1}$, $\tilde{v}_{2}$, and $S_{I_{4}}[w_1]$ in columns six to eight for the different parameter settings as well.
Moreover, the last columns of both tables present the cardinality of the full grids $\hat{G}_{N}^{21}$, where the sFFT algorithms search for the frequencies of the sparse representations of the computed approximations.

\begin{Example}
\begin{sloppy}
We consider the average error of the computed approximations $\breve{u}^\rho_\dagger$, $\dagger\in\{\rol,\mrol\}$, of the solution $u^*$ for fixed spatial nodes $\eta_k$.
To this end we calculate the solution of \eqref{eq:num_ode} for $n_{\mathrm{test}}=20000$ fixed randomly chosen parameters $\boldxi^{i} \in [-1,1]^{20},\, i=1, \dots, 20000$ as grid functions defined on the uniform grid
\begin{align}
\eta_{k}=\frac{k}{100}, \qquad k=0, \dots, 100,\label{eq:xgrid}
\end{align}
via numerical integration and an error bound of $10^{-6}$. We denote the corresponding solution by $\check{u}$, i.e., we assume that the values $\check{u}(\eta_k,\boldxi^i)$ are suitable approximations of the true solution and we use these function values for comparison against our approximations.  For a first comparison, we consider the pointwise difference with respect to our approximated solution and calculate the mean, i.e.
\begin{align*}
\Err^{\rho}_\dagger(\eta_{k}):=\frac{1}{n_{\mathrm{test}}}\sum_{i=1}^{n_{\mathrm{test}}}|\check{u}(\eta_{k},\boldxi^i)-\breve{u}_{\dagger}^{\rho}(\eta_{k},\boldxi^{i})|.
\end{align*}
The errors $\Err^{\rho}_\rol$ and $\Err^{\rho}_\mrol$ for the parameter selections $\rho=\rn{1}, \rn{2}, \rn{3}$ from Tables \ref{tab:parameter_and_nosn_sFFT} and \ref{tab:parameter_and_nosn_msFFT} are plotted in Figures \ref{fig:Errj_sFFT} and
\ref{fig:Errj_msFFT}, respectively.
We observe that the approximations computed by the MR1LsFFT are slightly better than those computed using the R1LsFFT algorithm. Certainly, this observation seems reasonable due to the usage of different numbers of sampling values, cf. Tables \ref{tab:parameter_and_nosn_msFFT} and \ref{tab:parameter_and_nosn_sFFT}.
\end{sloppy}

\begin{figure}[tb]
\begin{tikzpicture}
\begin{axis}[ymode=log,legend style={at={(0.05,0.05)},anchor=south west, /tikz/every even column/.append style={column sep=10pt}},legend columns = 3, width=1\textwidth,
       height=0.4\textwidth, xmin=0, xmax=1, ymin=5e-6, ymax=2e-3 ]
\addplot[mark=none,  black, dotted, very thick] table[x index=0,y index=1] {avrb.txt};
\addlegendentry{$\mathrm{Err}^{\rn{1}}_\rol$}
\addplot[mark=none, black, dashdotted,very thick] table[x index=0,y index=2] {avrb.txt};
\addlegendentry{$\mathrm{Err}^{\rn{2}}_\rol$}
\addplot[mark=none, black, dashed, very thick] table[x index=0,y index=3] {avrb.txt};
\addlegendentry{$\mathrm{Err}^{\rn{3}}_\rol$}
\end{axis}
\end{tikzpicture}
\caption{Averaged absolute errors $\Err^\rho_\rol$ for 20000 random samples of $y$ and using the sFFT algorithm with single rank-1 lattices sampling.}
\label{fig:Errj_sFFT}
\end{figure}

\begin{figure}[tb]
\begin{tikzpicture}
\begin{axis}[ymode=log,legend style={at={(0.05,0.05)},anchor=south west, /tikz/every even column/.append style={column sep=10pt}},legend columns = 3, width=1\textwidth,
       height=0.4\textwidth, xmin=0, xmax=1, ymin=5e-6, ymax=2e-3  ]
\addplot[mark=none,  black, dotted, very thick] table[x index=0,y index=1] {avrbmlft.txt};
\addlegendentry{$\mathrm{Err}^{\rn{1}}_\mrol$}
\addplot[mark=none,  black, dashdotted, very thick] table[x index=0,y index=2] {avrbmlft.txt};
\addlegendentry{$\mathrm{Err}^{\rn{2}}_\mrol$}
\addplot[mark=none,  black, dashed, very thick] table[x index=0,y index=3] {avrbmlft.txt};
\addlegendentry{$\mathrm{Err}^{\rn{3}}_\mrol$}
\end{axis}
\end{tikzpicture}
\caption{Averaged absolute errors $\Err^\rho_\mrol$ for 20000 random samples of $y$ and using the sFFT algorithm based on multiple rank-1 lattice sampling.}
\label{fig:Errj_msFFT}
\end{figure}
\end{Example}

According to the last example, we computed a complete approximate solution of \eqref{eq:num_ode}.
In Section \ref{sec:comp_moments} we explained how to compute moments of these approximate solution.

\begin{Example}
We demonstrate the performance of our approximation strategy by a comparison of subsequently computed approximate moments of the solution $u^{*}$. The Monte-Carlo approximation of the expectation value is given by
\begin{equation*}
\overline{u_{n_{\mathrm{test}}}(\eta_k)}=\frac {1}{n_{\mathrm{test}}} \sum_{i=1}^{n_{\mathrm{test}}} \check{u}(\eta_k,\boldxi^i)
\end{equation*}
for fixed $\eta_k$, cf. \eqref{eq:xgrid},
and the pointwise error at these spatial nodes is computed by
$$\Res^{\rho}_\dagger(\eta_{k}):=|\overline{u_{n_{\mathrm{test}}}(\eta_k)}-\mathbb{E}\breve{u}_{\dagger}^{\rho}(\eta_{k})|,$$
where the approximations $\mathbb{E}\breve{u}_{\dagger}^{\rho}(\eta_{k})$ of the first moment are gained from the solutions $\breve{u}_{\dagger}^{\rho}$, $\dagger\in\{\rol,\mrol\}$ as described in Section \ref{sec:comp_moments}.
The $\Res^{\rho}_\dagger$ behave very similar for fixed $\rho$ and $\dagger\in\{\rol,\mrol\}$. Slightly better errors of the Expectation can be observed for the multiple rank-1 lattice approach, cf. Figures \ref{fig:Resj_sFFT} and \ref{fig:Resj_msFFT}.

\begin{figure}[tb]
\begin{tikzpicture}
\begin{axis}[ymode=log,
legend style={at={(0.25,0.65)}, anchor=south west, /tikz/every even column/.append style={column sep=10pt}},legend columns = 4,
width=1\textwidth,
       height=0.4\textwidth, xmin=0, xmax=1, ymin=5e-7, ymax=1 ]
       \addplot[mark=none,  black, solid] coordinates{
			(0.00,1e-10) (0.01,1.165e-02) (0.02,2.307e-02) (0.03,3.424e-02) (0.04,4.518e-02) (0.05,5.589e-02) (0.06,6.635e-02) (0.07,7.658e-02) (0.08,8.658e-02) (0.09,9.634e-02) (0.10,1.059e-01) (0.11,1.151e-01) (0.12,1.242e-01) (0.13,1.330e-01) (0.14,1.416e-01) (0.15,1.499e-01) (0.16,1.580e-01) (0.17,1.659e-01) (0.18,1.735e-01) (0.19,1.809e-01) (0.20,1.881e-01) (0.21,1.950e-01) (0.22,2.017e-01) (0.23,2.082e-01) (0.24,2.144e-01) (0.25,2.204e-01) (0.26,2.261e-01) (0.27,2.316e-01) (0.28,2.369e-01) (0.29,2.420e-01) (0.30,2.468e-01) (0.31,2.514e-01) (0.32,2.557e-01) (0.33,2.598e-01) (0.34,2.637e-01) (0.35,2.673e-01) (0.36,2.707e-01) (0.37,2.739e-01) (0.38,2.768e-01) (0.39,2.795e-01) (0.40,2.819e-01) (0.41,2.842e-01) (0.42,2.862e-01) (0.43,2.879e-01) (0.44,2.894e-01) (0.45,2.907e-01) (0.46,2.918e-01) (0.47,2.926e-01) (0.48,2.932e-01) (0.49,2.935e-01) (0.50,2.937e-01) (0.51,2.935e-01) (0.52,2.932e-01) (0.53,2.926e-01) (0.54,2.918e-01) (0.55,2.907e-01) (0.56,2.894e-01) (0.57,2.879e-01) (0.58,2.862e-01) (0.59,2.842e-01) (0.60,2.819e-01) (0.61,2.795e-01) (0.62,2.768e-01) (0.63,2.739e-01) (0.64,2.707e-01) (0.65,2.673e-01) (0.66,2.637e-01) (0.67,2.598e-01) (0.68,2.557e-01) (0.69,2.514e-01) (0.70,2.468e-01) (0.71,2.420e-01) (0.72,2.369e-01) (0.73,2.316e-01) (0.74,2.261e-01) (0.75,2.204e-01) (0.76,2.144e-01) (0.77,2.082e-01) (0.78,2.017e-01) (0.79,1.951e-01) (0.80,1.881e-01) (0.81,1.810e-01) (0.82,1.736e-01) (0.83,1.659e-01) (0.84,1.581e-01) (0.85,1.500e-01) (0.86,1.416e-01) (0.87,1.331e-01) (0.88,1.242e-01) (0.89,1.152e-01) (0.90,1.059e-01) (0.91,9.637e-02) (0.92,8.661e-02) (0.93,7.662e-02) (0.94,6.638e-02) (0.95,5.591e-02) (0.96,4.520e-02) (0.97,3.426e-02) (0.98,2.308e-02) (0.99,1.166e-02) (1.00,1e-10)
              };
\addlegendentry{$\overline{u_{n_{\mathrm{test}}}(\eta_k)}$}
\addplot[mark=none,  black, dotted, very thick] table[x index=0,y index=1] {errewrb.txt};
\addlegendentry{$\Res^{\rn{1}}_\rol$}
\addplot[mark=none,  black, dashdotted, very thick] table[x index=0,y index=2] {errewrb.txt};
\addlegendentry{$\Res^{\rn{2}}_\rol$}
\addplot[mark=none,  black, dashed, very thick] table[x index=0,y index=3] {errewrb.txt};
\addlegendentry{$\Res^{\rn{3}}_\rol$}
\end{axis}
\end{tikzpicture}
\caption{Absolute difference $\Res^{\rho}_\rol(\eta_{k})$ of Monte-Carlo expectation value $\overline{u_{n_{\mathrm{test}}}(\eta_k)}$ and the approximation $\mathbb{E}\breve{u}^\rho_\rol$ of $u^*_{\operatorname{E}}$.}\label{fig:Resj_sFFT}
\end{figure}

\begin{figure}[hbt]
\begin{tikzpicture}
\begin{axis}[ymode=log,
legend style={at={(0.25,0.65)}, anchor=south west, /tikz/every even column/.append style={column sep=10pt}},legend columns = 4,
width=1\textwidth,
       height=0.4\textwidth, xmin=0, xmax=1, ymin=5e-7, ymax=1 ]
       \addplot[mark=none,  black, solid] coordinates{
			(0.00,1e-10) (0.01,1.165e-02) (0.02,2.307e-02) (0.03,3.424e-02) (0.04,4.518e-02) (0.05,5.589e-02) (0.06,6.635e-02) (0.07,7.658e-02) (0.08,8.658e-02) (0.09,9.634e-02) (0.10,1.059e-01) (0.11,1.151e-01) (0.12,1.242e-01) (0.13,1.330e-01) (0.14,1.416e-01) (0.15,1.499e-01) (0.16,1.580e-01) (0.17,1.659e-01) (0.18,1.735e-01) (0.19,1.809e-01) (0.20,1.881e-01) (0.21,1.950e-01) (0.22,2.017e-01) (0.23,2.082e-01) (0.24,2.144e-01) (0.25,2.204e-01) (0.26,2.261e-01) (0.27,2.316e-01) (0.28,2.369e-01) (0.29,2.420e-01) (0.30,2.468e-01) (0.31,2.514e-01) (0.32,2.557e-01) (0.33,2.598e-01) (0.34,2.637e-01) (0.35,2.673e-01) (0.36,2.707e-01) (0.37,2.739e-01) (0.38,2.768e-01) (0.39,2.795e-01) (0.40,2.819e-01) (0.41,2.842e-01) (0.42,2.862e-01) (0.43,2.879e-01) (0.44,2.894e-01) (0.45,2.907e-01) (0.46,2.918e-01) (0.47,2.926e-01) (0.48,2.932e-01) (0.49,2.935e-01) (0.50,2.937e-01) (0.51,2.935e-01) (0.52,2.932e-01) (0.53,2.926e-01) (0.54,2.918e-01) (0.55,2.907e-01) (0.56,2.894e-01) (0.57,2.879e-01) (0.58,2.862e-01) (0.59,2.842e-01) (0.60,2.819e-01) (0.61,2.795e-01) (0.62,2.768e-01) (0.63,2.739e-01) (0.64,2.707e-01) (0.65,2.673e-01) (0.66,2.637e-01) (0.67,2.598e-01) (0.68,2.557e-01) (0.69,2.514e-01) (0.70,2.468e-01) (0.71,2.420e-01) (0.72,2.369e-01) (0.73,2.316e-01) (0.74,2.261e-01) (0.75,2.204e-01) (0.76,2.144e-01) (0.77,2.082e-01) (0.78,2.017e-01) (0.79,1.951e-01) (0.80,1.881e-01) (0.81,1.810e-01) (0.82,1.736e-01) (0.83,1.659e-01) (0.84,1.581e-01) (0.85,1.500e-01) (0.86,1.416e-01) (0.87,1.331e-01) (0.88,1.242e-01) (0.89,1.152e-01) (0.90,1.059e-01) (0.91,9.637e-02) (0.92,8.661e-02) (0.93,7.662e-02) (0.94,6.638e-02) (0.95,5.591e-02) (0.96,4.520e-02) (0.97,3.426e-02) (0.98,2.308e-02) (0.99,1.166e-02) (1.00,1e-10)
              };
\addlegendentry{$\overline{u_{n_{\mathrm{test}}}(\eta_k)}$}
\addplot[mark=none, black, dotted, very thick] table[x index=0,y index=1] {errewrbmlft.txt};
\addlegendentry{$\Res^{\rn{1}}_\mrol$}
\addplot[mark=none, black, dashdotted, very thick] table[x index=0,y index=2] {errewrbmlft.txt};
\addlegendentry{$\Res^{\rn{2}}_\mrol$}
\addplot[mark=none, black, dashed, very thick] table[x index=0,y index=3] {errewrbmlft.txt};
\addlegendentry{$\Res^{\rn{3}}_\mrol$}
\end{axis}
\end{tikzpicture}
\caption{Absolute difference $\Res^{\rho}_\mrol(\eta_{k})$ of Monte-Carlo expectation value $\overline{u_{n_{\mathrm{test}}}(\eta_k)}$ and the approximation $\mathbb{E}\breve{u}^\rho_\mrol$ of $u^*_{\operatorname{E}}$.}\label{fig:Resj_msFFT}
\end{figure}

Furthermore, we can regard higher order moments. In a similar way as above, we computed the approximation of the second order moment by averaging
\begin{equation*}
\overline{u_{n_{\mathrm{test}}}(\eta_k)^2}=\frac {1}{n_{\mathrm{test}}} \sum_{i=1}^{n_{\mathrm{test}}} \check{u}(\eta_k,\boldxi^i)^2
\end{equation*}
and the pointwise error for each $\eta_k$
$$
\Res^{\rho,2}_\dagger(\eta_{k}):=|\overline{u_{n_{\mathrm{test}}}(\eta_k)^2}-\mathbb{E}\breve{u}_{\dagger}^{\rho}(\eta_{k})^2|.
$$

\begin{figure}[tb]
\begin{tikzpicture}
\begin{axis}[ymode=log,
legend style={at={(0.05,0.05)},anchor=south west, /tikz/every even column/.append style={column sep=10pt}},legend columns = 4,
width=1\textwidth,
       height=0.4\textwidth, xmin=0, xmax=1, ymin=1e-7, ymax=1e-1 ]
\addplot[mark=none,  black, solid] coordinates{
       (0.00e+00,1.093e-33) (1.00e-02,1.369e-04) (2.00e-02,5.363e-04) (3.00e-02,1.182e-03) (4.00e-02,2.058e-03) (5.00e-02,3.148e-03) (6.00e-02,4.437e-03) (7.00e-02,5.910e-03) (8.00e-02,7.553e-03) (9.00e-02,9.350e-03) (1.00e-01,1.129e-02) (1.10e-01,1.335e-02) (1.20e-01,1.553e-02) (1.30e-01,1.782e-02) (1.40e-01,2.018e-02) (1.50e-01,2.263e-02) (1.60e-01,2.514e-02) (1.70e-01,2.770e-02) (1.80e-01,3.031e-02) (1.90e-01,3.294e-02) (2.00e-01,3.560e-02) (2.10e-01,3.826e-02) (2.20e-01,4.093e-02) (2.30e-01,4.359e-02) (2.40e-01,4.622e-02) (2.50e-01,4.884e-02) (2.60e-01,5.141e-02) (2.70e-01,5.394e-02) (2.80e-01,5.642e-02) (2.90e-01,5.885e-02) (3.00e-01,6.120e-02) (3.10e-01,6.349e-02) (3.20e-01,6.569e-02) (3.30e-01,6.781e-02) (3.40e-01,6.984e-02) (3.50e-01,7.177e-02) (3.60e-01,7.361e-02) (3.70e-01,7.533e-02) (3.80e-01,7.695e-02) (3.90e-01,7.845e-02) (4.00e-01,7.983e-02) (4.10e-01,8.109e-02) (4.20e-01,8.223e-02) (4.30e-01,8.324e-02) (4.40e-01,8.412e-02) (4.50e-01,8.487e-02) (4.60e-01,8.548e-02) (4.70e-01,8.597e-02) (4.80e-01,8.631e-02) (4.90e-01,8.652e-02) (5.00e-01,8.659e-02) (5.10e-01,8.652e-02) (5.20e-01,8.631e-02) (5.30e-01,8.597e-02) (5.40e-01,8.549e-02) (5.50e-01,8.488e-02) (5.60e-01,8.413e-02) (5.70e-01,8.326e-02) (5.80e-01,8.225e-02) (5.90e-01,8.111e-02) (6.00e-01,7.985e-02) (6.10e-01,7.847e-02) (6.20e-01,7.697e-02) (6.30e-01,7.536e-02) (6.40e-01,7.363e-02) (6.50e-01,7.180e-02) (6.60e-01,6.987e-02) (6.70e-01,6.784e-02) (6.80e-01,6.572e-02) (6.90e-01,6.352e-02) (7.00e-01,6.123e-02) (7.10e-01,5.887e-02) (7.20e-01,5.645e-02) (7.30e-01,5.397e-02) (7.40e-01,5.144e-02) (7.50e-01,4.886e-02) (7.60e-01,4.625e-02) (7.70e-01,4.361e-02) (7.80e-01,4.095e-02) (7.90e-01,3.828e-02) (8.00e-01,3.562e-02) (8.10e-01,3.296e-02) (8.20e-01,3.032e-02) (8.30e-01,2.772e-02) (8.40e-01,2.515e-02) (8.50e-01,2.264e-02) (8.60e-01,2.020e-02) (8.70e-01,1.783e-02) (8.80e-01,1.554e-02) (8.90e-01,1.336e-02) (9.00e-01,1.130e-02) (9.10e-01,9.357e-03) (9.20e-01,7.558e-03) (9.30e-01,5.915e-03) (9.40e-01,4.441e-03) (9.50e-01,3.151e-03) (9.60e-01,2.060e-03) (9.70e-01,1.183e-03) (9.80e-01,5.369e-04) (9.90e-01,1.370e-04) (1.00e+00,1.093e-33)
       };
       \addlegendentry{$\overline{u_{n_{\mathrm{test}}}(\eta_k)^2}$}
\addplot[mark=none,  black, dotted, very thick] coordinates{
(0.000e+00,1.820e-06) (1.000e-02,1.492e-05) (2.000e-02,2.972e-05) (3.000e-02,4.923e-05) (4.000e-02,7.060e-05) (5.000e-02,8.975e-05) (6.000e-02,1.044e-04) (7.000e-02,1.146e-04) (8.000e-02,1.213e-04) (9.000e-02,1.239e-04) (1.000e-01,1.200e-04) (1.100e-01,1.070e-04) (1.200e-01,8.467e-05) (1.300e-01,5.753e-05) (1.400e-01,3.469e-05) (1.500e-01,2.654e-05) (1.600e-01,3.971e-05) (1.700e-01,7.284e-05) (1.800e-01,1.156e-04) (1.900e-01,1.521e-04) (2.000e-01,1.676e-04) (2.100e-01,1.553e-04) (2.200e-01,1.196e-04) (2.300e-01,7.469e-05) (2.400e-01,3.802e-05) (2.500e-01,2.224e-05) (2.600e-01,2.972e-05) (2.700e-01,5.199e-05) (2.800e-01,7.471e-05) (2.900e-01,8.548e-05) (3.000e-01,8.055e-05) (3.100e-01,6.680e-05) (3.200e-01,5.795e-05) (3.300e-01,6.697e-05) (3.400e-01,9.860e-05) (3.500e-01,1.460e-04) (3.600e-01,1.938e-04) (3.700e-01,2.251e-04) (3.800e-01,2.307e-04) (3.900e-01,2.136e-04) (4.000e-01,1.879e-04) (4.100e-01,1.717e-04) (4.200e-01,1.778e-04) (4.300e-01,2.067e-04) (4.400e-01,2.463e-04) (4.500e-01,2.776e-04) (4.600e-01,2.845e-04) (4.700e-01,2.624e-04) (4.800e-01,2.213e-04) (4.900e-01,1.817e-04) (5.000e-01,1.645e-04) (5.100e-01,1.817e-04) (5.200e-01,2.301e-04) (5.300e-01,2.928e-04) (5.400e-01,3.476e-04) (5.500e-01,3.769e-04) (5.600e-01,3.757e-04) (5.700e-01,3.531e-04) (5.800e-01,3.263e-04) (5.900e-01,3.116e-04) (6.000e-01,3.157e-04) (6.100e-01,3.322e-04) (6.200e-01,3.456e-04) (6.300e-01,3.393e-04) (6.400e-01,3.048e-04) (6.500e-01,2.460e-04) (6.600e-01,1.776e-04) (6.700e-01,1.181e-04) (6.800e-01,8.091e-05) (6.900e-01,6.830e-05) (7.000e-01,7.154e-05) (7.100e-01,7.620e-05) (7.200e-01,7.032e-05) (7.300e-01,5.096e-05) (7.400e-01,2.563e-05) (7.500e-01,8.216e-06) (7.600e-01,1.132e-05) (7.700e-01,3.926e-05) (7.800e-01,8.522e-05) (7.900e-01,1.339e-04) (8.000e-01,1.684e-04) (8.100e-01,1.777e-04) (8.200e-01,1.615e-04) (8.300e-01,1.293e-04) (8.400e-01,9.588e-05) (8.500e-01,7.411e-05) (8.600e-01,6.977e-05) (8.700e-01,8.024e-05) (8.800e-01,9.708e-05) (8.900e-01,1.109e-04) (9.000e-01,1.157e-04) (9.100e-01,1.107e-04) (9.200e-01,9.882e-05) (9.300e-01,8.427e-05) (9.400e-01,6.963e-05) (9.500e-01,5.528e-05) (9.600e-01,4.069e-05) (9.700e-01,2.629e-05) (9.800e-01,1.414e-05) (9.900e-01,6.161e-06) (1.000e+00,1.578e-07)
};
\addlegendentry{$\Res^{\rn{1},2}_\rol$}
\addplot[mark=none,  black, dashdotted, very thick] coordinates{
(0.000e+00,6.773e-08) (1.000e-02,2.139e-06) (2.000e-02,1.150e-06) (3.000e-02,8.155e-07) (4.000e-02,1.910e-06) (5.000e-02,3.885e-06) (6.000e-02,7.797e-06) (7.000e-02,1.164e-05) (8.000e-02,1.341e-05) (9.000e-02,1.406e-05) (1.000e-01,1.618e-05) (1.100e-01,2.075e-05) (1.200e-01,2.596e-05) (1.300e-01,2.867e-05) (1.400e-01,2.720e-05) (1.500e-01,2.367e-05) (1.600e-01,2.294e-05) (1.700e-01,2.735e-05) (1.800e-01,3.303e-05) (1.900e-01,3.414e-05) (2.000e-01,3.068e-05) (2.100e-01,2.869e-05) (2.200e-01,3.138e-05) (2.300e-01,3.392e-05) (2.400e-01,3.069e-05) (2.500e-01,2.439e-05) (2.600e-01,2.254e-05) (2.700e-01,2.634e-05) (2.800e-01,2.890e-05) (2.900e-01,2.596e-05) (3.000e-01,2.228e-05) (3.100e-01,2.407e-05) (3.200e-01,2.975e-05) (3.300e-01,3.321e-05) (3.400e-01,3.347e-05) (3.500e-01,3.507e-05) (3.600e-01,3.981e-05) (3.700e-01,4.411e-05) (3.800e-01,4.556e-05) (3.900e-01,4.649e-05) (4.000e-01,4.874e-05) (4.100e-01,4.955e-05) (4.200e-01,4.648e-05) (4.300e-01,4.251e-05) (4.400e-01,4.134e-05) (4.500e-01,4.032e-05) (4.600e-01,3.359e-05) (4.700e-01,2.167e-05) (4.800e-01,1.180e-05) (4.900e-01,7.807e-06) (5.000e-01,5.201e-06) (5.100e-01,1.047e-06) (5.200e-01,8.726e-06) (5.300e-01,1.247e-05) (5.400e-01,1.212e-05) (5.500e-01,1.208e-05) (5.600e-01,1.416e-05) (5.700e-01,1.631e-05) (5.800e-01,1.765e-05) (5.900e-01,1.960e-05) (6.000e-01,2.127e-05) (6.100e-01,1.934e-05) (6.200e-01,1.420e-05) (6.300e-01,1.156e-05) (6.400e-01,1.427e-05) (6.500e-01,1.697e-05) (6.600e-01,1.362e-05) (6.700e-01,7.307e-06) (6.800e-01,6.457e-06) (6.900e-01,1.231e-05) (7.000e-01,1.697e-05) (7.100e-01,1.519e-05) (7.200e-01,1.191e-05) (7.300e-01,1.416e-05) (7.400e-01,2.057e-05) (7.500e-01,2.414e-05) (7.600e-01,2.293e-05) (7.700e-01,2.178e-05) (7.800e-01,2.375e-05) (7.900e-01,2.572e-05) (8.000e-01,2.406e-05) (8.100e-01,2.044e-05) (8.200e-01,1.860e-05) (8.300e-01,1.834e-05) (8.400e-01,1.612e-05) (8.500e-01,1.055e-05) (8.600e-01,4.225e-06) (8.700e-01,2.804e-07) (8.800e-01,6.201e-07) (8.900e-01,1.185e-07) (9.000e-01,6.984e-07) (9.100e-01,2.981e-08) (9.200e-01,1.494e-06) (9.300e-01,2.284e-06) (9.400e-01,1.632e-06) (9.500e-01,9.439e-07) (9.600e-01,1.496e-06) (9.700e-01,1.950e-06) (9.800e-01,6.125e-07) (9.900e-01,6.859e-07) (1.000e+00,2.315e-09)
};
\addlegendentry{$\Res^{\rn{2},2}_\rol$}
\addplot[mark=none,  black, dashed, very thick] coordinates{
(0.000e+00,3.364e-08) (1.000e-02,5.362e-07) (2.000e-02,5.287e-07) (3.000e-02,8.275e-07) (4.000e-02,7.518e-07) (5.000e-02,2.631e-06) (6.000e-02,5.169e-06) (7.000e-02,6.053e-06) (8.000e-02,6.624e-06) (9.000e-02,8.852e-06) (1.000e-01,1.102e-05) (1.100e-01,1.080e-05) (1.200e-01,9.378e-06) (1.300e-01,8.739e-06) (1.400e-01,8.009e-06) (1.500e-01,6.111e-06) (1.600e-01,4.613e-06) (1.700e-01,4.584e-06) (1.800e-01,3.931e-06) (1.900e-01,1.133e-06) (2.000e-01,1.709e-06) (2.100e-01,2.612e-06) (2.200e-01,3.061e-06) (2.300e-01,4.576e-06) (2.400e-01,5.391e-06) (2.500e-01,3.772e-06) (2.600e-01,1.480e-06) (2.700e-01,8.638e-07) (2.800e-01,1.492e-06) (2.900e-01,2.419e-06) (3.000e-01,4.609e-06) (3.100e-01,8.261e-06) (3.200e-01,1.080e-05) (3.300e-01,1.037e-05) (3.400e-01,8.721e-06) (3.500e-01,8.079e-06) (3.600e-01,8.217e-06) (3.700e-01,8.521e-06) (3.800e-01,9.831e-06) (3.900e-01,1.192e-05) (4.000e-01,1.236e-05) (4.100e-01,1.016e-05) (4.200e-01,7.503e-06) (4.300e-01,6.043e-06) (4.400e-01,5.007e-06) (4.500e-01,4.526e-06) (4.600e-01,6.724e-06) (4.700e-01,1.133e-05) (4.800e-01,1.448e-05) (4.900e-01,1.450e-05) (5.000e-01,1.465e-05) (5.100e-01,1.746e-05) (5.200e-01,2.048e-05) (5.300e-01,2.099e-05) (5.400e-01,2.104e-05) (5.500e-01,2.393e-05) (5.600e-01,2.921e-05) (5.700e-01,3.450e-05) (5.800e-01,3.939e-05) (5.900e-01,4.402e-05) (6.000e-01,4.662e-05) (6.100e-01,4.588e-05) (6.200e-01,4.350e-05) (6.300e-01,4.180e-05) (6.400e-01,4.109e-05) (6.500e-01,4.141e-05) (6.600e-01,4.387e-05) (6.700e-01,4.796e-05) (6.800e-01,5.069e-05) (6.900e-01,5.076e-05) (7.000e-01,5.052e-05) (7.100e-01,5.175e-05) (7.200e-01,5.287e-05) (7.300e-01,5.254e-05) (7.400e-01,5.203e-05) (7.500e-01,5.182e-05) (7.600e-01,4.949e-05) (7.700e-01,4.386e-05) (7.800e-01,3.749e-05) (7.900e-01,3.256e-05) (8.000e-01,2.774e-05) (8.100e-01,2.178e-05) (8.200e-01,1.643e-05) (8.300e-01,1.313e-05) (8.400e-01,1.039e-05) (8.500e-01,6.965e-06) (8.600e-01,4.178e-06) (8.700e-01,2.682e-06) (8.800e-01,7.847e-07) (8.900e-01,1.774e-06) (9.000e-01,2.568e-06) (9.100e-01,1.163e-06) (9.200e-01,5.561e-07) (9.300e-01,1.763e-06) (9.400e-01,1.922e-06) (9.500e-01,3.190e-08) (9.600e-01,8.910e-07) (9.700e-01,1.920e-07) (9.800e-01,3.873e-07) (9.900e-01,5.305e-07) (1.000e+00,1.060e-09)
};
\addlegendentry{$\Res^{\rn{3},2}_\rol$}
\end{axis}
\end{tikzpicture}
\caption{Absolute difference $\Res^{\rho,2}_\rol(\eta_{k})$ of Monte-Carlo second moment $\overline{u_{n_{\mathrm{test}}}(\eta_k)^2}$ and the approximation $\mathbb{E}(\breve{u}^\rho_\rol)^2$ of $u^*_{\operatorname{E}^2}$.}\label{fig:Res2j_sFFT}
\end{figure}
In Figure \ref{fig:Res2j_sFFT}, the errors $\Res^{\rho,2}_\rol$ for the single rank-1 lattice approach (R1LsFFT)
are plotted. We see that adequately chosen parameters yield even very well approximated second moments.
\end{Example}

\subsection*{Acknowledgement}
LK and DP gratefully acknowledge the funding by the Deutsche Forschungsgemeinschaft (DFG, German Research Foundation, Projektnummer -- 380648269).


\small

\begin{thebibliography}{10}

\bibitem{BaCoDVGi17}
M.~Bachmayr, A.~Cohen, R.~DeVore, and G.~Migliorati.
\newblock Sparse polynomial approximation of parametric elliptic {PDEs}. part
  {II}: lognormal coefficients.
\newblock {\em ESAIM: M2AN}, 51:341--363, 2017.

\bibitem{BaCoGi17}
M.~Bachmayr, A.~Cohen, and G.~Migliorati.
\newblock Sparse polynomial approximation of parametric elliptic {PDEs}. part
  {I}: affine coefficients.
\newblock {\em ESAIM: M2AN}, 51:321--339, 2017.

\bibitem{BoRaSchw17}
J.-L. {Bouchot}, H.~{Rauhut}, and C.~{Schwab}.
\newblock {Multi-level Compressed Sensing Petrov-Galerkin discretization of
  high-dimensional parametric PDEs}.
\newblock {\em ArXiv e-prints}, Jan. 2017.
\newblock arXiv:1701.01671 [math.NA].

\bibitem{ChIwKr18}
B.~{Choi}, M.~{Iwen}, and F.~{Krahmer}.
\newblock Sparse harmonic transforms: A new class of sublinear-time algorithms
  for learning functions of many variables.
\newblock {\em ArXiv e-prints 1808.04932}, 2018.

\bibitem{CoDeSchw10a}
A.~Cohen, R.~DeVore, and C.~Schwab.
\newblock Analytic regularity and polynomial approximation of parametric and
  stochastic elliptic {PDE}s.
\newblock {\em Anal. Appl. (Singap.)}, 9:11 -- 47, 2010.

\bibitem{CoKuNuSu16}
R.~Cools, F.~Y. Kuo, D.~Nuyens, and G.~Suryanarayana.
\newblock {Tent-transformed lattice rules for integration and approximation of
  multivariate non-periodic functions}.
\newblock {\em J. Complexity}, 36:166--181, 2016.

\bibitem{EiPfSchn15}
M.~Eigel, M.~Pfeffer, and R.~Schneider.
\newblock Adaptive stochastic {G}alerkin {FEM} with hierarchical tensor
  representations.
\newblock {\em Numer. Math.}, 136:765--803, 2017.

\bibitem{GrKuNuScheSl11}
I.~Graham, F.~Kuo, D.~Nuyens, R.~Scheichl, and I.~Sloan.
\newblock Quasi-{M}onte {C}arlo methods for elliptic {PDE}s with random
  coefficients and applications.
\newblock {\em J. Comp. Phys.}, 230:3668 -- 3694, 2011.

\bibitem{HaNoTaTe16}
A.-L. Haji-Ali, F.~Nobile, L.~Tamellini, and R.~Tempone.
\newblock Multi-index stochastic collocation for random {PDEs}.
\newblock {\em Comput. Methods Appl. Mech. Engrg.}, 306:95--122, 2016.

\bibitem{HaSchw13}
M.~Hansen and C.~Schwab.
\newblock Analytic regularity and nonlinear approximation of a class of
  parametric semilinear elliptic {PDEs}.
\newblock {\em Math. Nachr.}, 286(8‐9):832--860, 2013.

\bibitem{kaemmererdiss}
L.~K{\"a}mmerer.
\newblock {\em {High Dimensional Fast Fourier Transform Based on Rank-1 Lattice
  Sampling}}.
\newblock Dissertation. Universit\"atsverlag Chemnitz, 2014.

\bibitem{Kae2013}
L.~K{\"a}mmerer.
\newblock Reconstructing multivariate trigonometric polynomials from samples
  along rank-1 lattices.
\newblock In G.~E. Fasshauer and L.~L. Schumaker, editors, {\em Approximation
  Theory XIV: San Antonio 2013}, pages 255--271. Springer International
  Publishing, 2014.

\bibitem{Kae17}
L.~{K{\"a}mmerer}.
\newblock {Constructing spatial discretizations for sparse multivariate
  trigonometric polynomials that allow for a fast discrete Fourier transform}.
\newblock {\em Appl. Comput. Harmon. Anal.}, 2017, accepted.

\bibitem{Kae16}
L.~K{\"{a}}mmerer.
\newblock Multiple rank-1 lattices as sampling schemes for multivariate
  trigonometric polynomials.
\newblock {\em J. Fourier Anal. Appl.}, 24:17--44, 2018.

\bibitem{KaPoVo17}
L.~{K{\"a}mmerer}, D.~{Potts}, and T.~{Volkmer}.
\newblock {High-dimensional sparse FFT based on sampling along multiple rank-1
  lattices}.
\newblock {\em ArXiv e-prints 1711.05152}, 2017.

\bibitem{MaKn10}
O.~P. {Le Ma\^{i}tre} and O.~M. Knio.
\newblock {\em Spectral Methods for Uncertainty Quantification}.
\newblock Scientific Computation. Springer Netherlands, 2010.

\bibitem{PoTa12}
D.~Potts and M.~Tasche.
\newblock Parameter estimation for multivariate exponential sums.
\newblock {\em Electron. Trans. Numer. Anal.}, 40:204--224, 2013.

\bibitem{PoVo14}
D.~Potts and T.~Volkmer.
\newblock {Sparse high-dimensional FFT based on rank-1 lattice sampling}.
\newblock {\em Appl. Comput. Harmon. Anal.}, 41:713--748, 2016.

\bibitem{RaSchw16}
H.~Rauhut and C.~Schwab.
\newblock Compressive sensing {P}etrov--{G}alerkin approximation of
  high-dimensional parametric operator equations.
\newblock {\em Math. Comp.}, 86:661--700, 2017.

\bibitem{SuNuCo14}
G.~Suryanarayana, D.~Nuyens, and R.~Cools.
\newblock Reconstruction and collocation of a class of non-periodic functions
  by sampling along tent-transformed rank-1 lattices.
\newblock {\em J. Fourier Anal. Appl.}, 22:187--214, 2016.

\bibitem{Vo_FFTr1l}
T.~Volkmer.
\newblock {sparseFFTr1l, \textsc{Matlab}$^{\text{\textregistered}}$ toolbox for
  computing the sparse fast Fourier transform based on reconstructing rank-1
  lattices in a dimension incremental way}.
\newblock \url{http://www.tu-chemnitz.de/~tovo/software}, 2015.

\end{thebibliography}
\end{document}